\documentclass{siamart1116}
\usepackage{amssymb}
\usepackage{eurosym}
\usepackage{lipsum}
\usepackage{amsfonts}
\usepackage{graphicx}
\usepackage{algorithmic}
\usepackage{amsopn}

\ifpdf
  \DeclareGraphicsExtensions{.eps,.pdf,.png,.jpg}
\else
  \DeclareGraphicsExtensions{.eps}
\fi
\numberwithin{theorem}{section}
\numberwithin{equation}{section}
\newcommand{\TheTitle}{Globally convergent numerical
method}
\newcommand{\TheAuthors}{M.V.Klibanov J.Li and W.Zhang}

\headers{\TheTitle}{\TheAuthors}

\ifpdf
\hypersetup{
  pdftitle={\TheTitle},
  pdfauthor={\TheAuthors}
}
\fi

\renewcommand{\thefigure}{\arabic{section}.\arabic{figure}}
\renewcommand\theequation{\thesection.\arabic{equation}}
\input{tcilatex}
\begin{document}

\title{A globally convergent numerical method for a 3D coefficient inverse
problem for a wave-like equation\thanks{%
Submitted to the editors DATE. 
\funding{
The work of Klibanov was supported by US Army Research Laboratory and US Army
Research Office grant W911NF-19-1-0044. The work of Li  was partially supported by the NSF of China No. 11971221 and the Shenzhen Sci-Tech Fund No. RCJC20200714114556020, JCYJ20170818153840322 and JCYJ20190809150413261,  
and Guangdong Provincial Key Laboratory of Computational Science
and Material Design No. 2019B030301001.
The work of Zhang was partially supported by the Shenzhen Sci-Tech Fund No. RCBS20200714114941241 and the NSF of China No.
11901282.}}}
\author{ Michael V. Klibanov\thanks{%
Department of Mathematics and Statistics, University of North Carolina at
Charlotte, Charlotte, NC 28223, USA (mklibanv@uncc.edu)} \and Jingzhi Li%
\thanks{%
Department of Mathematics, Southern University of Science and Technology
(SUSTech), 1088 Xueyuan Boulevard, University Town of Shenzhen, Xili,
Nanshan, Shenzhen, Guangdong Province, P.R.China (li.jz@sustech.edu.cn)}
\and Wenlong Zhang\thanks{%
Department of Mathematics, Southern University of Science and Technology
(SUSTech), 1088 Xueyuan Boulevard, University Town of Shenzhen, Xili,
Nanshan, Shenzhen, Guangdong Province, P.R.China (zhangwl@sustech.edu.cn)} }
\date{}
\maketitle

\begin{abstract}
A version of the convexification globally convergent numerical method is
constructed for a coefficient inverse problem for a wave-like partial
differential equation. The presence of the Carleman Weight Function in the
corresponding Tikhonov-like cost functional ensures the global strict
convexity of this functional. Numerical results are presented to illustrate
the effectiveness and efficiency of the proposed method.
\end{abstract}

\begin{keywords}
Carleman Weight Function, convexification, global
convergence, numerical results
\end{keywords}

\begin{AMS}
  35R30
\end{AMS}

\section{Introduction}

\label{sec:1}

In \cite{BKN} a version of the convexification globally convergent numerical
method was analytically developed for a 3-D Coefficient Inverse Problem
(CIP) for a wave-like Partial Differential Equation (PDE). In this paper we
first provide some new analytical results, which significantly enhance the
theory of \cite{BKN}. The main goal of this paper is to test the numerical
performance of the method of \cite{BKN}. A numerical study was not conducted
in \cite{BKN}. An application of the CIP of this paper to the elasticity
theory is discussed in the end of section 2.

Theorems 4.2 and 5.2 are new analytical results of this paper. The ideas of
proofs of these results are also new. Even though a result, similar with the
one of Theorem 5.2 was established in \cite[Theorem 4]{BKN}, we lift here
some restrictive conditions of \cite{BKN}. In Theorem 4.2, we prove the
Lipschitz stability estimate for a boundary value problem (BVP) for a
nonlinear PDE with a non local term. The presence of this term causes both
the nonlinearity and significant difficulties. Actually, we numerically
solve this BVP by the convexification method. Therefore, the importance of
the Lipschitz stability estimate is that it ensures that our BVP is
sufficiently stable. The latter is confirmed in numerical Tests 2 and 4 in
section 6, since these tests treat noisy data. 

The key step of the convexification method is the construction of a weighted
Tikhonov-like cost functional with the Carleman Weight Function (CWF) in it.
The CWF is involved as the weight function in the Carleman estimate for a
corresponding hyperbolic operator. This functional is strictly convex on a
convex bounded set $P$ in a Hilbert space. A smallness condition is not
imposed on the diameter of $P$, which, in turn means that our numerical
method converges globally. Unlike \cite{BKN}, we do not require here the
existence of the minimizer on $P$ of our functional. In addition, unlike 
\cite{BKN}, we do not assume here that this minimizer belongs to the
interior of $P$.

In addition, we do not impose anymore the assumption of Theorem 3 of \cite%
{BKN} that all terms of the sequence generated by the gradient method of the
minimization of that functional belong to the above mentioned bounded set.
Also, we decrease the required smoothness of the solution of the forward
problem from $H^{8}$ of \cite{BKN} to $H^{6}.$ In addition, the forward
problem we consider here is an initial boundary value problem for our PDE
rather than the Cauchy problem of \cite{BKN}. Finally, we present our
numerical results, which is our main goal here, and which was not done in 
\cite{BKN}.

\textbf{Remarks 1.1}:

\begin{enumerate}
\item \emph{It is well known that CIPs are both nonlinear and ill-posed.
These two factors cause substantial challenges in their studies. Because of
those challenges, as a rule, the minimal smoothness requirements are not of
a primary concern of the authors of corresponding publications, see, e.g. 
\cite{Nov2} as well as Theorem 4.1 in \cite{Rom}.}

\item \emph{In addition, our computational experience tells us that the
smoothness conditions can be significantly relaxed. This is related to both
the current paper and previous publications about the convexification cited
below.}
\end{enumerate}

The convexification method addresses the question of the globally convergent
numerical methods for CIPs with non overdetermined data. This is an
important question since conventional numerical methods for CIPs rely on the
minimizations of least squares cost functionals, see, e.g. \cite%
{Chavent,Gonch1,Gonch2}. However, these functionals are, as a rule, non
convex. The latter leads to the well known phenomenon of local minima and
ravines of those functionals. This phenomenon, in turn, means that in order
to get a good approximation for the exact solution of a CIP, one needs to
start the iterative optimization process of that functional in a
sufficiently small neighborhood of that solution. However, such a
neighborhood is rarely available in applications, also see \cite%
{Baud1,Baud2,Hoop} for similar comments. On the other hand, the global
convergence of the convexification method is rigorously guaranteed, see,
e.g. Theorem 5.3 as one of examples of such results as well as publications
about the convexification cited below.

We call a numerical method for a CIP \emph{globally convergent} if a theorem
is proven, which guarantees that this method reaches a sufficiently small
neighborhood of the true solution of that CIP without any advanced knowledge
of this neighborhood. Recall that one of the main concepts of the
regularization theory is the assumption of the existence of the true
solution with the \textquotedblleft ideal" noiseless data \cite{BK,T}. We
also refer to \cite{Hoop} for another globally convergent method for a
similar CIP.

The convexification numerical method was first proposed in 1995 \cite%
{KlibIous} and then in 1997 \cite{Klib97,KlibNW}. Some computations were
performed in \cite{KT} in the 1D case. Active computational studies of the
convexification method have started after the publication of the paper \cite%
{Bak} in 2017, since this work has analytically addressed a number of
questions about the numerics. We refer to \cite%
{Khoa1,Khoa2,KLZ1,KLZ2,KLNSN,KlibSAR2} and references therein for some
samples of publications, in which various analytical results about the
convexification are combined with numerical studies, including the cases of
experimentally collected backscattering data. We also refer to the recently
published book of the first two authors \cite{KL}. All these works consider
the cases when the wave field is generated by a point source. Even though
these publications are generated by the ideas of the Bukhgeim-Klibanov
method (BK), the case of the point source is outside of the framework of BK.
We refer to \cite{BukhKlib} for the originating paper on BK as well as to,
e.g. \cite{BK,BY,Ksurvey,KL} and references cited therein for some follow up
publications. In \cite{BukhKlib}, the method of Carleman estimates was
introduced in the field of CIPs for the first time. The focus of \cite%
{BukhKlib} was only on the question of uniqueness theorems for
multidimensional CIPs. First extensions of the idea of BK to the numerical
side of the theory of CIPs were published later \cite{KlibIous,Klib97,KlibNW}%
.

The framework of BK works only with the case when one of initial conditions
in a wave-like PDE does not equal zero in the entire domain of interest.
This assumption is used in \cite{BKN}. Therefore, we also use it in the
current paper, since we study here the numerical performance of the method
of \cite{BKN}. We also refer to works \cite{Baud1,Baud2}, where a different
version of the convexification method is used for CIPs for wave-like PDEs.
Both publications \cite{Baud1,Baud2} work within the framework of the BK\
method.

The paper is structured as follows. In section 2 we formulate our
Coefficient Inverse Problem. In section 3 we derive a BVP for a nonlinear
PDE with a non-local term. If this problem is solved, then the target
unknown coefficient can be easily found. In section 4 we formulate a
Carleman estimate and prove the above mentioned Lipschitz stability estimate
for that BVP. In section 5 we introduce the central functional of the
convexification method and investigate its properties. Section 6 is devoted
to a description of our numerical implementation as well as to the
demonstration of numerical results.

\section{Statement of the Coefficient Inverse Problem}

\label{sec:2}

In all Hilbert spaces considered below functions are real valued ones. Let $%
\Omega \subset \mathbb{R}^{3}$ be a convex bounded domain with a piecewise
smooth boundary $\partial \Omega \in C^{\infty }$ and let the number $T>0.$
Denote $Q_{T}=\Omega \times (0,T),S_{T}=\partial \Omega \times \left(
0,T\right) .$ Let the function $c\left( x\right) $ satisfies the following
conditions: 
\begin{equation}
c\left( x\right) \in \left[ 1,b\right] \text{ for }x\in \Omega ,  \label{2.1}
\end{equation}%
\begin{equation}
c\in C^{1}\left( \overline{\Omega }\right) .  \label{2.2}
\end{equation}%
We use \textquotedblleft 1" here in (\ref{2.1}) for the normalization only.
Here $b>0$ is \ fixed number. In addition to (\ref{2.1}), (\ref{2.2}), we
assume that there exists a point $x_{0}\in \mathbb{R}^{3}\diagdown \overline{%
\Omega }$ such that%
\begin{equation}
\min_{\overline{\Omega }}\left( \nabla c,x-x_{0}\right) \geq 0,  \label{4.1}
\end{equation}%
where $\left( ,\right) $ denotes the scalar product in $\mathbb{R}^{3}.$
Condition (\ref{4.1}) means that the function $c\left( x\right) $ is
increasing along the ray connecting any point $x\in \overline{\Omega }$ with
the point $x_{0}$ when $x$ moves away from $x_{0}.$ This condition is the
standard one in Carleman estimates for the hyperbolic operator $c\left(
x\right) \partial _{t}^{2}-\Delta ,$ see \cite[Theorem 1.10.2]{BK}, \cite[%
Theorem 2.5.1]{KL}. It is unclear whether or not that Carleman estimate is
valid without condition (\ref{4.1}).

We assume that 
\begin{equation}
f\in C^{3}\left( \overline{\Omega }\right) .  \label{2.06}
\end{equation}%
Let $a,k>0,k>a$ be two numbers. We assume below that 
\begin{equation}
k\geq \Delta f\geq a=const.>0\text{ for }x\in \overline{\Omega }\text{.}
\label{2.6}
\end{equation}%
Consider the following initial boundary value problem for a wave-like PDE:%
\begin{equation}
c\left( x\right) u_{tt}=\Delta u,\left( x,t\right) \in Q_{T},  \label{2.7}
\end{equation}%
\begin{equation}
u\left( x,0\right) =f\left( x\right) ,u_{t}\left( x,0\right) =0,  \label{2.8}
\end{equation}%
\begin{equation}
u|_{S_{T}}=s(x,t).  \label{2.9}
\end{equation}%
It is shown in \cite[Chapter 4, Corollary 4.1]{Lad} that one can impose some
non restrictive conditions on $\partial \Omega $ and functions $f\left(
x\right) $ and $s(x,t),$ which guarantee existence and uniqueness of the
solution 
\begin{equation}
u\in H^{6}\left( Q_{T}\right)  \label{2.10}
\end{equation}%
of problem (\ref{2.7})-(\ref{2.9}). However, we do not discuss these
conditions here for brevity. Rather, we assume smoothness (\ref{2.10}) of
the solution of that problem, also, see Remark 1.1 in section 1.

\textbf{Coefficient Inverse Problem (CIP). }\emph{Suppose that conditions (%
\ref{2.1})-(\ref{2.10}) are satisfied. Assume that the coefficient }$c(x)$%
\emph{\ is unknown inside of the domain }$\Omega $\emph{. Determine the
function }$c(x)$\emph{\ for }$x\in \Omega ,$\emph{\ assuming that the
Neumann boundary condition }$p\left( x,t\right) $\emph{\ is known, }%
\begin{equation}
\partial _{n}u\mid _{S_{T}}=p\left( x,t\right) ,  \label{2.11}
\end{equation}%
\emph{where }$n$\emph{\ is the normal outward looking vector at }$S_{T}.$

Uniqueness of this CIP \ under conditions (\ref{2.6}), (\ref{2.10}) was
established by the method of \cite{BukhKlib}, see, e.g. \cite[Theorem
1.10.5.1]{BK} and \cite[Theorem 3.2.1, Theorem 3.3]{KL}.

Equation (\ref{2.7}) is the acoustic equation, and the function $1/\sqrt{%
c\left( x\right) }$ is the speed of propagation of sound waves, and $u\left(
x,t\right) $ is the amplitude of sound waves. Equation (\ref{2.7}) also
appears in the elasticity theory for isotropic media \cite{Lurie}. Consider
a simplified case when $u\left( x,t\right) $ is one of components of the
displacement field, other components approximately equal zero and 
\begin{equation}
c\left( x\right) =\frac{\rho \left( x\right) }{\lambda \left( x\right) +2\mu
\left( x\right) },  \label{2.12}
\end{equation}%
where $\rho \left( x\right) $ is the density of the medium, $\lambda \left(
x\right) $ and $\mu \left( x\right) $ are Lam\'{e} coefficients \cite{Lurie}%
. Let $s(x,t)$ be the displacement at the boundary $\partial \Omega $ of the
medium and $p\left( x,t\right) $ be the normal stress at that boundary.
Suppose that we have a 3-D elastic slab occupying the domain $\Omega .$
Suppose that, applying a certain force at the moment of time $\left\{
t=0\right\} $, we arrange the shape of this slab to become $u\left(
x,0\right) =f\left( x\right) ,$ and the slab does not fluctuate at $\left\{
t=0\right\} ,$ i.e. $u_{t}\left( x,0\right) =0$. For times $t\in \left(
0,T\right) $ this slab is allowed to fluctuate without applying a force
inside of $\Omega $, although the displacement at the boundary of the slab
is controlled by the given function $s\left( x,t\right) =u\mid _{S_{T}}.$
Then we obtain forward problem (\ref{2.7})-(\ref{2.9}). Assume that we
measure the normal stress at the boundary of $\Omega $, i.e. we measure the
function $p\left( x,t\right) =\partial _{n}u\mid _{S_{T}}$. Then CIP (\ref%
{2.7})-(\ref{2.9}), (\ref{2.11}) is the problem of determining combination (%
\ref{2.12}) of the above named parameters using boundary measurements (\ref%
{2.11}).

\section{Nonlinear Equation With a Non-Local Term}

\label{sec:3}

In this section we derive from (\ref{2.7})-(\ref{2.9}), (\ref{2.11}) a
nonlinear BVP with a non local term. This BVP does not contain the unknown
coefficient $c\left( x\right) .$ Having the solution of this BVP, one can
straightforwardly compute $c\left( x\right) .$ Therefore, we focus on this
BVP below.

\subsection{Nonlinear boundary value problem}

\label{sec:3.1}

Denote 
\begin{equation}
q\left( x,t\right) =\partial _{t}^{2}s\left( x,t\right) ,r\left( x,t\right)
=\partial _{t}^{2}p\left( x,t\right) .  \label{3.1}
\end{equation}%
Let $w=u_{tt}.$ By (\ref{2.10}) $w\in H^{4}\left( Q_{T}\right) .$ Since by (%
\ref{2.7}) 
\begin{equation}
w\left( x,0\right) =\frac{\left( \Delta f\right) \left( x\right) }{c\left(
x\right) },  \label{3.2}
\end{equation}%
then (\ref{2.1}) and (\ref{2.6}) imply%
\begin{equation}
\frac{a}{b}\leq w\left( x,0\right) \leq k,  \label{3.3}
\end{equation}%
\begin{equation}
\frac{a}{k}\leq \frac{\Delta f}{w\left( x,0\right) }\leq \frac{kb}{a}.
\label{3.4}
\end{equation}%
We now establish an upper estimate on the constant $k.$ Fix an arbitrary
number $R>0$. We will seek for the function $w$ satisfying%
\begin{equation}
w\in H^{4}\left( Q_{T}\right) ,\text{ }\left\Vert w\right\Vert _{H^{4}\left(
Q_{T}\right) }\leq R.  \label{3.5}
\end{equation}%
By the embedding theorem 
\begin{equation}
H^{4}\left( Q_{T}\right) \subset C^{1}\left( \overline{Q}_{T}\right) \text{
and }\left\Vert v\right\Vert _{C^{1}\left( \overline{Q}_{T}\right) }\leq
C\left\Vert v\right\Vert _{H^{4}\left( Q_{T}\right) },\text{ }\forall v\in
H^{4}\left( Q_{T}\right) .  \label{3.50}
\end{equation}%
Hence, for all functions $w$ satisfying (\ref{3.5}) 
\begin{equation}
\left\Vert w\right\Vert _{C^{1}\left( \overline{Q}_{T}\right) }\leq CR,
\label{3.6}
\end{equation}%
where the number $C=C\left( Q_{T}\right) >0$ in (\ref{3.50}), (\ref{3.6})
depends only on the domain $Q_{T}.$ Hence, (\ref{3.3}) implies that we
should have $k\in \left( 0,CR\right] .$ Thus, we replace (\ref{3.3}), (\ref%
{3.4}) with: 
\begin{equation}
\frac{a}{CR}\leq \frac{\left( \Delta f\right) \left( x\right) }{w\left(
x,0\right) }\leq CR\frac{b}{a}.  \label{3.8}
\end{equation}%
Using (\ref{2.7})-(\ref{2.9}), (\ref{2.11}) and (\ref{3.1}) we obtain 
\begin{equation}
c\left( x\right) =\frac{\left( \Delta f\right) \left( x\right) }{w\left(
x,0\right) },\text{ }x\in \Omega ,  \label{3.9}
\end{equation}%
\begin{equation}
\frac{\left( \Delta f\right) \left( x\right) }{w\left( x,0\right) }%
w_{tt}-\Delta w=0,\text{ }\left( x,t\right) \in Q_{T},  \label{3.10}
\end{equation}%
\begin{equation}
w_{t}\left( x,0\right) =0,  \label{3.11}
\end{equation}%
\begin{equation}
w\mid _{S_{T}}=q\left( x,t\right) ,\text{ }\partial _{n}w\mid
_{S_{T}}=r\left( x,t\right) .  \label{3.12}
\end{equation}

Thus, $w\left( x,0\right) $ is the non local term in the PDE (\ref{3.10}).
In addition, PDE (\ref{3.10}) is nonlinear due to the presence of this term.
We focus below on the numerical solution of BVP (\ref{3.10})-(\ref{3.12}).

\subsection{The set for $w\left( x,t\right) $}

\label{sec:3.2}

Define the spaces $H_{0}^{2}\left( Q_{T}\right) ,H_{0}^{4}\left(
Q_{T}\right) ,H_{0}^{4,0}\left( Q_{T}\right) $ and the set $\Omega _{T}$ as: 
\begin{equation}
H_{0}^{2}\left( Q_{T}\right) =\left\{ u\in H^{2}\left( Q_{T}\right)
:u_{t}\left( x,0\right) =0\right\} ,  \label{4.4}
\end{equation}%
\[
H_{0}^{4}\left( Q_{T}\right) =\left\{ u\in H^{4}\left( Q_{T}\right)
:u_{t}\left( x,0\right) =0\right\} , 
\]%
\begin{equation}
H_{0}^{4,0}\left( Q_{T}\right) =\left\{ u\in H^{2}\left( Q_{T}\right)
:u_{t}\left( x,0\right) =0,u\mid _{S_{T}}=0,\partial _{n}u\mid
_{S_{T}}=0\right\} ,  \label{4.40}
\end{equation}%
\[
\Omega _{T}=\left\{ \left( x,t\right) :x\in \Omega ,t=T\right\} . 
\]%
We now specify the set of functions on which we search for the solution $w$
of BVP (\ref{3.10})-(\ref{3.12}). Denote%
\begin{equation}
A\left( w\right) \left( x\right) =\frac{\Delta f\left( x\right) }{w\left(
x,0\right) }.  \label{3.120}
\end{equation}%
By (\ref{2.06}), (\ref{3.5}) and (\ref{3.50}) the function $A\left( w\right)
\left( x\right) \in C^{1}\left( \overline{\Omega }\right) .$ We search for
the function $w\in P=P\left( a,b,R,f,q,r\right) \subset H_{0}^{4}\left(
Q_{T}\right) ,$ where the set $P$ depends only on listed parameters and is
defined as:%
\begin{equation}
P=P\left( x_{0},a,b,R,f,q,r\right) =\left\{ 
\begin{array}{c}
w\left( x,t\right) :w\in H_{0}^{4}\left( Q_{T}\right) , \\ 
\left\Vert w\right\Vert _{H^{4}\left( Q_{T}\right) }\leq R, \\ 
a/CR\leq A\left( w\right) \leq CRb/a, \\ 
\min_{\overline{\Omega }}\left( \nabla A\left( w\right) \left( x\right)
,x-x_{0}\right) \geq 0, \\ 
w\mid _{S_{T}}=q\left( x,t\right) ,\text{ }\partial _{n}w\mid
_{S_{T}}=r\left( x,t\right) .%
\end{array}%
\right.  \label{3.13}
\end{equation}%
The inequality in the third line of (\ref{3.13}) is due to (\ref{3.8}) and (%
\ref{3.120}). The inequality in the fourth line of (\ref{3.13}) follows from
(\ref{2.2}), (\ref{4.1}), (\ref{3.5}), (\ref{3.50}), (\ref{3.9}) and (\ref%
{3.120}). Obviously 
\begin{equation}
P=\overline{P}.  \label{3.130}
\end{equation}

\textbf{Lemma 3.1 }\cite{BKN}. \emph{The set }$P=P\left( a,b,R,f,q,r\right) $%
\emph{\ is convex. }

Uniqueness of the solution $w\in P$ of problem (\ref{3.10})-(\ref{3.12})
follows from Theorem 4.2. Suppose that we have computed the solution $w_{%
\text{comp}}\in \overline{P}$ of this problem. Then, using (\ref{3.2}) and (%
\ref{3.120}), we set the computed coefficient $c_{\text{comp}}\left(
x\right) $ as:%
\begin{equation}
c_{\text{comp}}\left( x\right) =\frac{\Delta f\left( x\right) }{w_{\text{comp%
}}\left( x,0\right) }.  \label{3.14}
\end{equation}

\section{The Carleman Estimate and the Lipschitz Stability Estimate}

\label{sec:4}

Since the convexification method is based on Carleman estimates, we
formulate in this section a Carleman estimate, which was proven in \cite{BKN}%
. Next, we formulate and prove the Lipschitz stability estimate for BVP (\ref%
{3.10})-(\ref{3.12}).

\subsection{Carleman estimate}

\label{sec:4.1}

Let the numbers $\eta \in \left( 0,1\right) $ and $\lambda \geq 1$. Consider
functions $\psi ,\varphi _{\lambda }$,%
\begin{equation}
\psi \left( x,t\right) =\left\vert x-x_{0}\right\vert ^{2}-\eta t^{2},\text{ 
}\varphi _{\lambda }\left( x,t\right) =\exp \left( \lambda \psi \left(
x,t\right) \right) .  \label{4.0}
\end{equation}%
Consider the numbers $d,D$,%
\begin{equation}
d=\min_{x\in \overline{\Omega }}\left\vert x-x_{0}\right\vert ,\text{ }%
D=\max_{x\in \overline{\Omega }}\left\vert x-x_{0}\right\vert .  \label{4.2}
\end{equation}%
For any number $\eta \in \left( 0,1\right) $ we choose a sufficiently large $%
T>D/\sqrt{\eta }.$ Then%
\begin{equation}
M=M\left( \Omega ,x_{0},\eta ,T\right) =\eta T^{2}-d^{2}>0,  \label{4.03}
\end{equation}%
\begin{equation}
N=N\left( \Omega ,x_{0},\eta ,T\right) =\eta T^{2}-D^{2}>0,  \label{4.3}
\end{equation}%
\begin{equation}
\varphi _{\lambda }^{2}\left( x,t\right) \leq e^{2\lambda D^{2}},\text{ }%
\left( x,t\right) \in Q_{T}.  \label{4.30}
\end{equation}%
Let the function $g\left( x,t\right) $ be such that $g,\partial _{t}g\in
C\left( \overline{Q}_{T}\right) .$

\textbf{Theorem 4.1 (}Carleman estimate \cite{BKN}). \emph{Assume that the
function }$c\left( x\right) $\emph{\ satisfies conditions (\ref{2.1})-(\ref%
{4.1}).\ Let }$D^{2}$ \emph{be the number defined in (\ref{4.2}).} \emph{%
Then there exist a number }$\eta _{0}=\eta _{0}\left( \Omega
,x_{0},\left\Vert c\right\Vert _{C^{1}\left( \overline{\Omega }\right)
}\right) \in \left( 0,1\right) $\emph{, a number }

$C_{1}=C_{1}\left( \eta ,x_{0},a,b,\left\Vert g\right\Vert _{C\left( 
\overline{Q}_{T}\right) },\left\Vert \partial _{t}g\right\Vert _{C\left( 
\overline{Q}_{T}\right) },Q_{T}\right) >0$\emph{\ and a sufficiently large
number }$\lambda _{0}=\lambda _{0}\left( \eta ,x_{0},a,b,\left\Vert
g\right\Vert _{C\left( \overline{Q}_{T}\right) },\left\Vert \partial
_{t}g\right\Vert _{C\left( \overline{Q}_{T}\right) },Q_{T}\right) \geq 1$%
\emph{, all three numbers depending only on listed parameters, such that if }%
$\eta \in \left( 0,\eta _{0}\right) $ \emph{and the number} $T=T\left( \eta
,x_{0},\Omega \right) >0$ \emph{is so large that (\ref{4.3}) holds, then for
all }$\lambda \geq \lambda _{0}$\emph{\ the following Carleman estimate is
valid}%
\[
\int\limits_{Q_{T}}\left( c\left( x\right) u_{tt}-\Delta u\right)
^{2}\varphi _{\lambda }^{2}dxdt+\int\limits_{Q_{T}}g\left( x,t\right)
u\left( x,0\right) u_{tt}\left( x,t\right) \varphi _{\lambda }^{2}dxdt 
\]%
\[
+C_{1}\lambda ^{3}\exp \left( 2\lambda D^{2}\right) \left( \left\Vert
u\right\Vert _{H^{1}\left( S_{T}\right) }^{2}+\left\Vert \partial _{n}u\mid
_{S_{T}}\right\Vert _{L_{2}\left( S_{T}\right) }^{2}\right) 
\]%
\[
+C_{1}\lambda ^{3}\exp \left( -2\lambda N\right) \left( \left\Vert
u_{t}\right\Vert _{L_{2}\left( \Omega _{T}\right) }^{2}+\left\Vert
u\right\Vert _{H^{1}\left( \Omega _{T}\right) }^{2}\right) +C_{1}\exp \left(
-2\lambda N\right) \left\Vert u\right\Vert _{H^{1}\left( Q_{T}\right) }^{2} 
\]%
\begin{equation}
\geq C_{1}\int\limits_{Q_{T}}\left( \lambda u_{t}^{2}+\lambda \left( \nabla
u\right) ^{2}+\lambda ^{3}u^{2}\right) \varphi _{\lambda }^{2}dxdt,\text{ }%
\forall u\in H_{0}^{2}\left( Q_{T}\right) .  \label{4.5}
\end{equation}

\textbf{Remarks 4.1:}

\begin{enumerate}
\item \emph{The nonlinear term }$g\left( x,t\right) u\left( x,0\right)
u_{tt}\left( x,t\right) $\emph{\ is introduced in (\ref{4.5}) because its
presence was used in \cite{BKN} the proof of a direct analog of Theorem 5.1
about the strict convexity of the functional, which is constructed in
section 5. If }$g\left( x,t\right) \equiv 0,$\emph{\ then (\ref{4.5}) is
just the conventional Carleman estimate of \cite[Theorem 1.10.2]{BK}, \cite[%
Theorem 2.5.1]{KL}. }

\item \emph{A non-conventional element here is the absence in (\ref{4.5}) of
the integral over }$\overline{Q}_{T}\cap \left\{ t=0\right\} .$\emph{\
Indeed, it follows from formulae (1.86), (1.87) of \cite{BK} as well as from
Corollary 2.5.1 of \cite{KL} that this integral equals to zero as long as }$%
u_{t}\left( x,0\right) =0,$\emph{\ see (\ref{4.4}).}

\item \emph{Since the function }$\varphi _{\lambda }^{2}\left( x,t\right) $%
\emph{\ is involved as the weight function in (\ref{4.5}), then we call }$%
\varphi _{\lambda }^{2}$\emph{\ the \textquotedblleft Carleman Weight
Function" (CWF) for the operator }$c\left( x\right) \partial _{t}^{2}-\Delta
.$
\end{enumerate}

\subsection{Lipschitz stability estimate}

\label{sec:4.2}

In this section, we use Theorem 4.1 to establish the Lipschitz stability
estimate for the nonlinear problem (\ref{3.10})-(\ref{3.12}). The
nonlinearity is due to the presence of the non local term $w\left(
x,0\right) $ in (\ref{3.10}). Suppose that we have two functions $%
w_{1},w_{2}\in H^{2}\left( Q_{T}\right) $ satisfying the following
conditions for $i=1,2:$ 
\begin{equation}
A\left( w_{i}\right) \left( x\right) \partial _{t}^{2}w_{i}-\Delta w_{i}=0,%
\text{ }\left( x,t\right) \in Q_{T},  \label{4.6}
\end{equation}%
\begin{equation}
\partial _{t}w_{i}\left( x,0\right) =0,i=1,2,  \label{4.7}
\end{equation}%
\begin{equation}
w_{i}\mid _{S_{T}}=q_{i}\left( x,t\right) ,\text{ }\partial _{n}w_{i}\mid
_{S_{T}}=r_{i}\left( x,t\right) ,  \label{4.8}
\end{equation}%
\begin{equation}
A\left( w_{i}\right) \left( x\right) =\frac{\Delta f\left( x\right) }{%
w_{i}\left( x,0\right) }.  \label{4.9}
\end{equation}

Recall that the space $H_{0}^{2}\left( Q_{T}\right) $ is defined in (\ref%
{4.4}).

\textbf{Theorem 4.2}. \emph{Assume that the function }$f$\emph{\ satisfies
condition (\ref{2.06}) and }

$\left\Vert f\right\Vert _{C^{3}\left( \overline{\Omega }\right) }\leq
f^{0}. $\emph{\ Assume that the function }$\Delta f\left( x\right) $\emph{\
satisfies conditions (\ref{2.6}) with two positive constants }$a,k,$\emph{\
where }$0<a<k.$\emph{For }$i=1,2$\emph{\ let the functions }$w_{i}\in
H_{0}^{2}\left( Q_{T}\right) $ \emph{satisfy conditions} \emph{(\ref{4.6})-(%
\ref{4.9}),} $w_{i}\left( x,0\right) \in C^{1}\left( \overline{\Omega }%
\right) $\emph{\ and }$\left\Vert w_{i}\left( x,0\right) \right\Vert
_{C^{1}\left( \overline{\Omega }\right) }\leq w^{0}.$\emph{\ Also, assume
that there exists a number }$m>0$\emph{\ such that }$w_{i}\left( x,0\right)
\geq m$\emph{\ in }$\overline{\Omega }.$\emph{\ In addition, suppose that
there exists a point }$x_{0}\in \mathbb{R}^{3}\diagdown \overline{\Omega }$%
\emph{\ such that }%
\begin{equation}
\left( \nabla A\left( w_{1}\right) \left( x\right) ,x-x_{0}\right) \geq 0,%
\text{ }\forall x\in \overline{\Omega }.  \label{4.10}
\end{equation}%
\emph{Let }$d$\emph{\ and }$D$\emph{\ be two numbers defined in (\ref{4.2})
and let the number }$T>D/\sqrt{\eta },$\emph{\ where the number }$\eta \in
\left( 0,\eta _{1}\right) $\emph{, where the number }$\eta _{1}=\eta
_{1}\left( \Omega ,x_{0},f^{0},w^{0},m\right) \in \left( 0,1\right) $ \emph{%
depends only on listed parameters and is chosen the same way as the number }$%
\eta _{0}\in \left( 0,1\right) $\emph{\ in Theorem 4.1, in which }$c\left(
x\right) $\emph{\ is replaced with }$A\left( w_{1}\right) \left( x\right) .$%
\emph{\ Then there exists a number }$Z=Z\left( \eta
,f^{0},w^{0},m,a,k,x_{0},Q_{T}\right) >0$\emph{\ depending only on listed
parameters such that the following} \emph{Lipschitz stability estimate holds:%
}%
\begin{equation}
\left\Vert w_{1}-w_{2}\right\Vert _{H^{1}\left( Q_{T}\right) }\leq Z\left(
\left\Vert q_{1}-q_{2}\right\Vert _{H^{1}\left( S_{T}\right) }+\left\Vert
r_{1}-r_{2}\right\Vert _{L_{2}\left( S_{T}\right) }\right) .  \label{4.11}
\end{equation}

\textbf{Proof}. Below in this proof $Z>0$ denotes different numbers
depending on parameters listed in the formulation of this theorem. It
follows from (\ref{3.4}) and (\ref{4.9}) that $A\left( w_{1}\right) \in
C^{1}\left( \overline{\Omega }\right) $ and 
\begin{equation}
\left\Vert A\left( w_{1}\right) \left( x\right) \right\Vert _{C^{1}\left( 
\overline{\Omega }\right) }\leq Z_{1},A\left( w_{1}\right) \left( x\right)
\geq \frac{a}{w^{0}}>0,  \label{4.110}
\end{equation}%
where the constant $Z_{1}=Z_{1}\left( f^{0},w^{0},m\right) >0$ depends only
on listed parameters. Denote 
\begin{equation}
\widetilde{w}=w_{1}-w_{2},\widetilde{q}=q_{1}-q_{2},\widetilde{r}%
=r_{1}-r_{2}.  \label{4.12}
\end{equation}%
Obviously, 
\begin{equation}
\widetilde{w}\left( x,0\right) =\widetilde{w}\left( x,t\right) -%
\mathop{\displaystyle \int}\limits_{0}^{t}\widetilde{w}_{t}\left( x,\tau
\right) d\tau .  \label{4.13}
\end{equation}%
Using the formula $a_{1}b_{1}-a_{2}b_{2}=\left( a_{1}-a_{2}\right)
b_{1}+\left( b_{1}-b_{2}\right) a_{2},\forall a_{1},b_{1},a_{2},b_{2}\in 
\mathbb{R},$ we obtain from (\ref{4.6})-(\ref{4.9}), (\ref{4.12}) and (\ref%
{4.13})

\begin{equation}
A\left( w_{1}\right) \left( x\right) \widetilde{w}_{tt}-\Delta \widetilde{w}%
+Y\left( x,t\right) \left( \widetilde{w}\left( x,t\right) -%
\mathop{\displaystyle \int}\limits_{0}^{t}\widetilde{w}_{t}\left( x,\tau
\right) d\tau \right) =0,  \label{4.14}
\end{equation}%
\begin{equation}
\widetilde{w}_{t}\left( x,0\right) =0,  \label{4.15}
\end{equation}%
\begin{equation}
\widetilde{w}\mid _{S_{T}}=\widetilde{q}\left( x,t\right) ,\text{ }\partial
_{n}\widetilde{w}\mid _{S_{T}}=\widetilde{r}\left( x,t\right) ,  \label{4.16}
\end{equation}%
\[
Y\left( x,t\right) =-\frac{\Delta f\left( x\right) }{w_{1}\left( x,t\right)
w_{2}\left( x,t\right) }. 
\]%
Hence, 
\begin{equation}
\left\vert Y^{2}\left( x,t\right) \right\vert \leq Z\text{ in }Q_{T}.
\label{4.17}
\end{equation}%
Let $\varphi _{\lambda }\left( x,t\right) ,\lambda \geq 1$ be the function
defined in (\ref{4.0}), where the number $\eta \in \left( 0,\eta _{1}\right) 
$\emph{\ }and the number $\eta _{1}=\eta _{1}\left( \Omega
,x_{0},f^{0},w^{0},m\right) \in \left( 0,1\right) $ is chosen the same way
as the number $\eta _{0}$ in Theorem 4.1, in which $c\left( x\right) $ is
replaced with $A\left( w_{1}\right) \left( x\right) .$ Using either Lemma
1.10.3 of \cite{BK} or Lemma 3.1.1 of \cite{KL} as well as (\ref{4.17}), we
obtain%
\[
\mathop{\displaystyle \int}\limits_{Q_{T}}Y^{2}\left( 
\mathop{\displaystyle
\int}\limits_{0}^{t}\widetilde{w}_{t}\left( x,\tau \right) d\tau \right)
^{2}\varphi _{\lambda }^{2}dxdt\leq \frac{Z}{\lambda }%
\mathop{\displaystyle
\int}\limits_{Q_{T}}\widetilde{w}_{t}^{2}\varphi _{\lambda }^{2}dxdt. 
\]%
Hence, (\ref{4.14}) and (\ref{4.17}) imply:%
\begin{equation}
\mathop{\displaystyle \int}\limits_{Q_{T}}\left( A\left( w_{1}\right) \left(
x\right) \widetilde{w}_{tt}-\Delta \widetilde{w}\right) ^{2}\varphi
_{\lambda }^{2}dxdt\leq Z\mathop{\displaystyle \int}\limits_{Q_{T}}\left( 
\widetilde{w}_{t}^{2}+\widetilde{w}^{2}\right) \varphi _{\lambda }^{2}dxdt.
\label{4.18}
\end{equation}

In Theorem 4.1, set $g\left( x,t\right) \equiv 0,$ see item 1 in Remarks
4.1. Then, keeping in mind (\ref{4.10}), (\ref{4.110}), (\ref{4.15}) and (%
\ref{4.16}), apply (\ref{4.5}) to the left hand side of (\ref{4.18}). We
obtain 
\[
Z\mathop{\displaystyle \int}\limits_{Q_{T}}\left( \widetilde{w}_{t}^{2}+%
\widetilde{w}^{2}\right) \varphi _{\lambda }^{2}dxdt+Z\lambda ^{3}\exp
\left( 2\lambda D^{2}\right) \left( \left\Vert \widetilde{q}\right\Vert
_{H^{1}\left( S_{T}\right) }^{2}+\left\Vert \widetilde{r}\right\Vert
_{L_{2}\left( S_{T}\right) }^{2}\right) 
\]%
\[
+Z\lambda ^{3}\exp \left( -2\lambda N\right) \left( \left\Vert \widetilde{w}%
_{t}\right\Vert _{L_{2}\left( \Omega _{T}\right) }^{2}+\left\Vert \widetilde{%
w}\right\Vert _{H^{1}\left( \Omega _{T}\right) }^{2}\right) +Z\exp \left(
-2\lambda N\right) \left\Vert \widetilde{w}\right\Vert _{H^{1}\left(
Q_{T}\right) }^{2} 
\]%
\begin{equation}
\geq \lambda \int\limits_{Q_{T}}\left( \widetilde{w}_{t}^{2}+\left( \nabla 
\widetilde{w}\right) ^{2}+\lambda ^{2}\widetilde{w}^{2}\right) \varphi
_{\lambda }^{2}dxdt,\text{ }\forall \lambda \geq \lambda _{0}\geq 1.
\label{4.19}
\end{equation}%
Choose the number $\lambda ^{\left( 1\right) }=\lambda ^{\left( 1\right)
}\left( \eta ,f^{0},w^{0},m,a,k,x_{0},Q_{T}\right) \geq \lambda _{0}\geq 1$
depending only on listed parameters such that $\lambda ^{\left( 1\right)
}\geq 2Z.$ Then the first term in the first line of (\ref{4.19}) is absorbed
by terms in the third line of (\ref{4.19}). Hence, 
\[
Z\exp \left( -3\lambda N/2\right) \left( \left\Vert \widetilde{w}%
_{t}\right\Vert _{L_{2}\left( \Omega _{T}\right) }^{2}+\left\Vert \widetilde{%
w}\right\Vert _{H^{1}\left( \Omega _{T}\right) }^{2}+\left\Vert \widetilde{w}%
\right\Vert _{H^{1}\left( Q_{T}\right) }^{2}\right) 
\]%
\begin{equation}
+Z\exp \left( 3\lambda D^{2}\right) \left( \left\Vert \widetilde{q}%
\right\Vert _{H^{1}\left( S_{T}\right) }^{2}+\left\Vert \widetilde{r}%
\right\Vert _{L_{2}\left( S_{T}\right) }^{2}\right)  \label{4.20}
\end{equation}%
\[
\geq \lambda \int\limits_{Q_{T}}\left( \widetilde{w}_{t}^{2}+\left( \nabla 
\widetilde{w}\right) ^{2}+\lambda ^{2}\widetilde{w}^{2}\right) \varphi
_{\lambda }^{2}dxdt,\text{ }\forall \lambda \geq \lambda ^{\left( 1\right)
}\geq 1. 
\]%
Choose the number 
\begin{equation}
t_{1}=\frac{1}{\sqrt{2}}\sqrt{T^{2}-D^{2}/\eta }<T.  \label{4.21}
\end{equation}%
Obviously $Q_{t_{1}}\subset Q_{T}$ and $\varphi _{\lambda }^{2}\left(
x,t\right) \geq -2\lambda \eta t_{1}^{2}$ for $\left( x,t\right) \in
Q_{t_{1}}.$ Also, using (\ref{4.3}) and (\ref{4.21}), we obtain%
\begin{equation}
-\frac{3}{2}\lambda N+2\lambda \eta t_{1}^{2}=-\frac{\lambda }{2}\left( \eta
T^{2}-D^{2}\right) =-\lambda \frac{N}{2}.  \label{4.22}
\end{equation}%
Hence, (\ref{4.20}) and (\ref{4.22}) imply%
\begin{equation}
\left\Vert \widetilde{w}\right\Vert _{H^{1}\left( Q_{t_{1}}\right) }^{2}\leq
Z\exp \left( -\lambda \frac{N}{2}\right) \left( \left\Vert \widetilde{w}%
_{t}\right\Vert _{L_{2}\left( \Omega _{T}\right) }^{2}+\left\Vert \widetilde{%
w}\right\Vert _{H^{1}\left( \Omega _{T}\right) }^{2}+\left\Vert \widetilde{w}%
\right\Vert _{H^{1}\left( Q_{T}\right) }^{2}\right)  \label{4.23}
\end{equation}%
\[
+Z\exp \left( 3\lambda \left( D^{2}+T^{2}\right) \right) \left( \left\Vert 
\widetilde{q}\right\Vert _{H^{1}\left( S_{T}\right) }^{2}+\left\Vert 
\widetilde{r}\right\Vert _{L_{2}\left( S_{T}\right) }^{2}\right) ,\forall
\lambda \geq \lambda ^{\left( 1\right) }\geq 1. 
\]%
Denote 
\begin{equation}
W\left( \lambda \right) =\text{the right hand side of (\ref{4.23}).}
\label{4.230}
\end{equation}%
By the mean value theorem there exists a number $t_{0}\in \left(
0,t_{1}\right) $ such that 
\begin{equation}
\left\Vert \widetilde{w}\left( x,t_{0}\right) \right\Vert _{H^{1}\left(
\Omega \right) }^{2}+\left\Vert \widetilde{w}_{t}\left( x,t_{0}\right)
\right\Vert _{L_{2}\left( \Omega \right) }^{2}=\frac{1}{t_{1}}\left\Vert 
\widetilde{w}\right\Vert _{H^{1}\left( Q_{t_{1}}\right) }^{2}\leq
Z\left\Vert \widetilde{w}\right\Vert _{H^{1}\left( Q_{t_{1}}\right) }^{2}.
\label{4.231}
\end{equation}%
Hence, by (\ref{4.23}) and (\ref{4.230}) 
\begin{equation}
\left\Vert \widetilde{w}\left( x,t_{0}\right) \right\Vert _{H^{1}\left(
\Omega \right) }^{2}+\left\Vert \widetilde{w}_{t}\left( x,t_{0}\right)
\right\Vert _{L_{2}\left( \Omega \right) }^{2}\leq W\left( \lambda \right) .
\label{4.24}
\end{equation}

We now apply the method of energy estimates \cite{Lad}\emph{\ }to problem (%
\ref{4.14})-(\ref{4.16}) in the domain $Q_{t_{0},T}=\Omega \times \left(
t_{0},T\right) .$ Let $x=\left( x_{1},x_{2},x_{3}\right) \in \mathbb{R}^{3}.$
Multiply both parts of equation (\ref{4.14}) by $2\widetilde{w}_{t}.$ Taking
into account (\ref{4.17}) and using Cauchy-Schwarz inequality, we obtain%
\[
\partial _{t}\left( A\left( w_{1}\right) \left( x\right) \widetilde{w}%
_{t}^{2}+\left( \nabla \widetilde{w}\right) ^{2}\right) +\mathop{%
\displaystyle \sum }\limits_{i=1}^{3}\partial _{x_{i}}\left( -2\widetilde{w}%
_{x_{i}}\widetilde{w}_{t}\right) 
\]%
\[
\leq Z\left( \widetilde{w}_{t}^{2}+\widetilde{w}^{2}+%
\mathop{\displaystyle
\int}\limits_{t_{0}}^{t}\widetilde{w}_{t}^{2}\left( x,\tau \right) d\tau +%
\mathop{\displaystyle \int}\limits_{0}^{t_{0}}\widetilde{w}_{t}^{2}\left(
x,\tau \right) d\tau \right) . 
\]%
Integrating this inequality over the domain $Q\left( t_{0},t^{\prime
}\right) $ for an arbitrary $t^{\prime }\in \left( t_{0},T\right) $ and
using the second inequality (\ref{4.110}), we obtain 
\begin{equation}
\mathop{\displaystyle \int}\limits_{\Omega }\left( \widetilde{w}%
_{t}^{2}+\left( \nabla \widetilde{w}\right) ^{2}\right) \left( x,t^{\prime
}\right) dx\leq Z\mathop{\displaystyle \int}\limits_{t_{0}}^{t^{\prime }}%
\mathop{\displaystyle \int}\limits_{\Omega }\left( \widetilde{w}_{t}^{2}+%
\widetilde{w}^{2}\right) \left( x,t\right) dxdt+Z\left\Vert \widetilde{w}%
_{t}\right\Vert _{L_{2}\left( Q_{t_{1}}\right) }^{2}  \label{4.25}
\end{equation}%
\[
+\left\Vert \widetilde{w}\left( x,t_{0}\right) \right\Vert _{H^{1}\left(
\Omega \right) }^{2}+\left\Vert \widetilde{w}_{t}\left( x,t_{0}\right)
\right\Vert _{L_{2}\left( \Omega \right) }^{2}+\left\Vert \widetilde{q}%
\right\Vert _{H^{1}\left( S_{T}\right) }^{2}+\left\Vert \widetilde{r}%
\right\Vert _{L_{2}\left( S_{T}\right) }^{2}+\left\Vert \widetilde{w}%
\right\Vert _{H^{1}\left( Q_{t_{1}}\right) }^{2}. 
\]%
Using (\ref{4.25}), Gronwall inequality and substituting in (\ref{4.13}) $%
t=t_{0}$ instead of $t=0,$we obtain%
\[
\mathop{\displaystyle \int}\limits_{\Omega }\left( \widetilde{w}%
_{t}^{2}+\left( \nabla \widetilde{w}\right) ^{2}+\widetilde{w}^{2}\right)
\left( x,t^{\prime }\right) dx\leq Z\left( \left\Vert \widetilde{w}\left(
x,t_{0}\right) \right\Vert _{H^{1}\left( \Omega \right) }^{2}+\left\Vert 
\widetilde{w}_{t}\left( x,t_{0}\right) \right\Vert _{L_{2}\left( \Omega
\right) }^{2}\right) 
\]%
\begin{equation}
+Z\left( \left\Vert \widetilde{w}\right\Vert _{H^{1}\left( Q_{t_{1}}\right)
}^{2}+\left\Vert \widetilde{q}\right\Vert _{H^{1}\left( S_{T}\right)
}^{2}+\left\Vert \widetilde{r}\right\Vert _{L_{2}\left( S_{T}\right)
}^{2}\right) ,\forall t^{\prime }\in \left( t_{0},T\right] .  \label{4.26}
\end{equation}%
Setting $t^{\prime }=T$ in (\ref{4.26}), we obtain%
\[
\left\Vert \widetilde{w}_{t}\right\Vert _{L_{2}\left( \Omega _{T}\right)
}^{2}+\left\Vert \widetilde{w}\right\Vert _{H^{1}\left( \Omega _{T}\right)
}^{2}\leq Z\left( \left\Vert \widetilde{q}\right\Vert _{H^{1}\left(
S_{T}\right) }^{2}+\left\Vert \widetilde{r}\right\Vert _{L_{2}\left(
S_{T}\right) }^{2}\right) 
\]%
\begin{equation}
+Z\left( \left\Vert \widetilde{w}\left( x,t_{0}\right) \right\Vert
_{H^{1}\left( \Omega \right) }^{2}+\left\Vert \widetilde{w}_{t}\left(
x,t_{0}\right) \right\Vert _{L_{2}\left( \Omega \right) }^{2}+\left\Vert 
\widetilde{w}\right\Vert _{H^{1}\left( Q_{t_{1}}\right) }^{2}\right) .
\label{4.27}
\end{equation}%
In addition, integrating (\ref{4.26}) with respect to $t^{\prime }\in \left(
t_{0},T\right] ,$ we obtain%
\begin{equation}
\left\Vert \widetilde{w}\right\Vert _{H^{1}\left( Q_{t_{0},T}\right)
}^{2}\leq Z\left( \left\Vert \widetilde{w}\left( x,t_{0}\right) \right\Vert
_{H^{1}\left( \Omega \right) }^{2}+\left\Vert \widetilde{w}_{t}\left(
x,t_{0}\right) \right\Vert _{L_{2}\left( \Omega \right) }^{2}+\left\Vert 
\widetilde{w}\right\Vert _{H^{1}\left( Q_{t_{1}}\right) }^{2}\right)
\label{4.28}
\end{equation}%
\[
+Z\left( \left\Vert \widetilde{q}\right\Vert _{H^{1}\left( S_{T}\right)
}^{2}+\left\Vert \widetilde{r}\right\Vert _{L_{2}\left( S_{T}\right)
}^{2}\right) . 
\]%
Substituting (\ref{4.230}) and (\ref{4.24}) in (\ref{4.27}) and (\ref{4.28}%
), we obtain%
\begin{equation}
\left\Vert \widetilde{w}_{t}\right\Vert _{L_{2}\left( \Omega _{T}\right)
}^{2}+\left\Vert \widetilde{w}\right\Vert _{H^{1}\left( \Omega _{T}\right)
}^{2}\leq W\left( \lambda \right) ,\text{ }\forall \lambda \geq \lambda
^{\left( 1\right) }\geq 1,  \label{4.29}
\end{equation}%
and also%
\begin{equation}
\left\Vert \widetilde{w}\right\Vert _{H^{1}\left( Q_{t_{0},T}\right)
}^{2}\leq W\left( \lambda \right) ,\text{ }\forall \lambda \geq \lambda
^{\left( 1\right) }\geq 1.  \label{4.31}
\end{equation}%
Let the number $\lambda ^{\left( 2\right) }=\lambda ^{\left( 2\right)
}\left( \eta ,f^{0},w^{0},m,a,k,x_{0},Q_{T}\right) \geq \lambda ^{\left(
1\right) }\geq 1$ be so large that 
\begin{equation}
Z\exp \left( -\lambda ^{\left( 2\right) }\frac{N}{2}\right) \leq \frac{1}{4}.
\label{4.310}
\end{equation}
Hence, (\ref{4.230}) and (\ref{4.29}) imply%
\[
\frac{3}{4}\left( \left\Vert \widetilde{w}_{t}\right\Vert _{L_{2}\left(
\Omega _{T}\right) }^{2}+\left\Vert \widetilde{w}\right\Vert _{H^{1}\left(
\Omega _{T}\right) }^{2}\right) \leq 
\]%
\[
\leq \left( 1-Z\exp \left( -\lambda ^{\left( 2\right) }\frac{N}{2}\right)
\right) \left( \left\Vert \widetilde{w}_{t}\right\Vert _{L_{2}\left( \Omega
_{T}\right) }^{2}+\left\Vert \widetilde{w}\right\Vert _{H^{1}\left( \Omega
_{T}\right) }^{2}\right) 
\]%
\[
\leq Z\exp \left( -\lambda ^{\left( 2\right) }\frac{N}{2}\right) \left\Vert 
\widetilde{w}\right\Vert _{H^{1}\left( Q_{T}\right) }^{2} 
\]%
\[
+Z\exp \left( 3\lambda ^{\left( 2\right) }\left( D^{2}+T^{2}\right) \right)
\left( \left\Vert \widetilde{q}\right\Vert _{H^{1}\left( S_{T}\right)
}^{2}+\left\Vert \widetilde{r}\right\Vert _{L_{2}\left( S_{T}\right)
}^{2}\right) . 
\]%
Hence, 
\begin{equation}
\left\Vert \widetilde{w}_{t}\right\Vert _{L_{2}\left( \Omega _{T}\right)
}^{2}+\left\Vert \widetilde{w}\right\Vert _{H^{1}\left( \Omega _{T}\right)
}^{2}\leq Z\exp \left( -\lambda ^{\left( 2\right) }\frac{N}{2}\right)
\left\Vert \widetilde{w}\right\Vert _{H^{1}\left( Q_{T}\right) }^{2}
\label{4.32}
\end{equation}%
\[
+Z\exp \left( 3\lambda ^{\left( 2\right) }\left( D^{2}+T^{2}\right) \right)
\left( \left\Vert \widetilde{q}\right\Vert _{H^{1}\left( S_{T}\right)
}^{2}+\left\Vert \widetilde{r}\right\Vert _{L_{2}\left( S_{T}\right)
}^{2}\right) . 
\]%
Using (\ref{4.23}), (\ref{4.230}), (\ref{4.31}) and (\ref{4.32}), we obtain%
\begin{equation}
\left\Vert \widetilde{w}\right\Vert _{H^{1}\left( Q_{t_{0},T}\right)
}^{2}\leq Z\exp \left( -\lambda ^{\left( 2\right) }\frac{N}{2}\right)
\left\Vert \widetilde{w}\right\Vert _{H^{1}\left( Q_{T}\right) }^{2}
\label{4.33}
\end{equation}%
\[
+Z\exp \left( 3\lambda ^{\left( 2\right) }\left( D^{2}+T^{2}\right) \right)
\left( \left\Vert \widetilde{q}\right\Vert _{H^{1}\left( S_{T}\right)
}^{2}+\left\Vert \widetilde{r}\right\Vert _{L_{2}\left( S_{T}\right)
}^{2}\right) . 
\]%
Furthermore, combining (\ref{4.23}) and (\ref{4.32}) and recalling that
since $t_{0}\in \left( 0,t_{1}\right) ,$ then $Q_{t_{0}}\subset Q_{t_{1}},$
we obtain%
\begin{equation}
\left\Vert \widetilde{w}\right\Vert _{H^{1}\left( Q_{t_{0}}\right) }^{2}\leq
\left\Vert \widetilde{w}\right\Vert _{H^{1}\left( Q_{t_{1}}\right) }^{2}\leq
Z\exp \left( -\lambda ^{\left( 2\right) }\frac{N}{2}\right) \left\Vert 
\widetilde{w}\right\Vert _{H^{1}\left( Q_{T}\right) }^{2}  \label{4.34}
\end{equation}%
\[
+Z\exp \left( 3\lambda ^{\left( 2\right) }\left( D^{2}+T^{2}\right) \right)
\left( \left\Vert \widetilde{q}\right\Vert _{H^{1}\left( S_{T}\right)
}^{2}+\left\Vert \widetilde{r}\right\Vert _{L_{2}\left( S_{T}\right)
}^{2}\right) . 
\]%
Obviously $\left\Vert \widetilde{w}\right\Vert _{H^{1}\left(
Q_{t_{0},T}\right) }^{2}+\left\Vert \widetilde{w}\right\Vert _{H^{1}\left(
Q_{t_{0}}\right) }^{2}=\left\Vert \widetilde{w}\right\Vert _{H^{1}\left(
Q_{T}\right) }^{2}.$ Hence, summing up (\ref{4.33}) and (\ref{4.34}) and
using (\ref{4.310}), we obtain%
\begin{equation}
\frac{1}{2}\left\Vert \widetilde{w}\right\Vert _{H^{1}\left( Q_{T}\right)
}^{2}\leq \left( 1-2Z\exp \left( -\lambda ^{\left( 2\right) }\frac{N}{2}%
\right) \right) \left\Vert \widetilde{w}\right\Vert _{H^{1}\left(
Q_{T}\right) }^{2}  \label{4.35}
\end{equation}%
\[
+Z\exp \left( 3\lambda ^{\left( 2\right) }\left( D^{2}+T^{2}\right) \right)
\left( \left\Vert \widetilde{q}\right\Vert _{H^{1}\left( S_{T}\right)
}^{2}+\left\Vert \widetilde{r}\right\Vert _{L_{2}\left( S_{T}\right)
}^{2}\right) . 
\]%
The target estimate (\ref{4.11}) of this theorem follows immediately from (%
\ref{4.12}) and (\ref{4.35}). $\square $

\section{A Globally Strictly Convex\ Cost Functional}

\label{sec:5}

In this section we present analytical results of our paper. To solve BVP (%
\ref{3.10})-(\ref{3.12}) numerically, we construct in this section a
Tikhonov-like functional, which is strictly convex on the set $P$ in (\ref%
{3.13}). Since a smallness condition is not imposed on the number $R$ in (%
\ref{3.13}), then this functional is \emph{globally} strictly convex.

\subsection{Construction of a globally strictly convex cost functional}

\label{sec:5.1}

Let $\alpha \in \left( 0,1\right) $ be the regularization parameter.
Recalling (\ref{3.130}),\emph{\ }we construct the following weighted
Tikhonov-like regularization functional $J_{\lambda ,\alpha }\left( w\right)
:P\rightarrow \mathbb{R}$: 
\begin{equation}
J_{\lambda ,\alpha }\left( w\right) =\int\limits_{Q_{T}}\left[ A\left(
w\right) \left( w\right) _{tt}-\Delta \left( w\right) \right] ^{2}\varphi
_{\lambda }^{2}dxdt+\alpha \left\Vert w\right\Vert _{H^{4}\left(
Q_{T}\right) }^{2},\text{ }w\in P.  \label{5.1}
\end{equation}

\textbf{Minimization Problem}. \emph{Minimize the functional }$J_{\lambda
,\alpha }\left( w\right) $\emph{\ on the set }$P=P\left( a,b,R,f,q,r\right) $%
\emph{\ defined in (\ref{3.13}).}

\textbf{Theorem 5.1}:

\emph{1. For any values of parameters }$\lambda ,\alpha ,\eta ,T>0$\emph{\
and for any function }$w\in P$\emph{\ there exists the Fr\'{e}chet
derivative }$J_{\lambda ,\alpha }^{\prime }\left( w\right) \in
H_{0}^{4,0}\left( Q_{T}\right) $\emph{\ of the functional }$J_{\lambda
,\alpha }.$\emph{\ Furthermore, this derivative is Lipschitz continuous,
i.e. there exists a constant }$K>0$\emph{\ such that }%
\begin{equation}
\left\Vert J_{\lambda ,\alpha }^{\prime }\left( w_{1}\right) -J_{\lambda
,\alpha }^{\prime }\left( w_{2}\right) \right\Vert _{H^{4}\left(
Q_{T}\right) }\leq K\left\Vert w_{1}-w_{2}\right\Vert _{H^{4}\left(
Q_{T}\right) },\forall w_{1},w_{2}\in H_{0}^{4}\left( Q_{T}\right) .
\label{5.2}
\end{equation}

\emph{2. Let }$\eta $\emph{,}$T$\emph{\ and }$\lambda _{0}$ \emph{be the
numbers of Theorem 4.1. Then there exists a constant }$C_{2}=C_{2}\left(
\eta ,x_{0},a,b,f,R,P,Q_{T}\right) >0$\emph{\ and a sufficiently large
number }$\lambda _{1}=\lambda _{1}\left( \eta ,x_{0},a,b,f,R,P,Q_{T}\right)
\geq \lambda _{0},$\emph{\ both depending only on listed parameters, such
that for all }$\lambda \geq \lambda _{1}$\emph{\ and for all }$\alpha \in %
\left[ 2e^{-\lambda N},1\right) $\emph{\ the functional }$J_{\lambda ,\alpha
}\left( w\right) $\emph{\ is strictly convex on the set }$P.$\emph{\ More
precisely,} \emph{\ }%
\[
J_{\lambda ,\alpha }\left( w_{2}\right) -J_{\lambda ,\alpha }\left(
w_{1}\right) -J_{\lambda ,\alpha }^{\prime }\left( w_{1}\right) \left(
w_{2}-w_{1}\right) 
\]%
\begin{equation}
\geq C_{2}e^{-2\lambda M}\left\Vert w_{2}-w_{1}\right\Vert _{H^{1}\left(
Q_{T}\right) }^{2}+\frac{\alpha }{2}\left\Vert \left( w_{2}-w_{1}\right)
\right\Vert _{H^{4}\left( Q_{T}\right) }^{2}  \label{5.3}
\end{equation}%
\[
\geq \frac{\alpha }{2}\left\Vert w_{2}-w_{1}\right\Vert _{H^{4}\left(
Q_{T}\right) }^{2},\text{ }\forall w_{1},w_{2}\in P, 
\]%
\emph{where the constant }$M>0$\emph{\ is defined in (\ref{4.03}).}

\emph{3. There exists unique minimizer }$w_{\min }\in P$\emph{\ of the
functional }$J_{\lambda ,\alpha }\left( w\right) $\emph{\ on the set }$P$%
\emph{\ and the following inequality holds:}%
\begin{equation}
J_{\lambda ,\alpha }^{\prime }\left( w_{\min }\right) \left( w-w_{\min
}\right) \geq 0,\text{ }\forall w\in P.  \label{5.4}
\end{equation}

Item 1 of this theorem is proven in \cite[Theorem 2]{BKN}, except of (\ref%
{5.2}). However,\emph{\ }(\ref{5.2}) can be proven completely similarly with
the proof of Theorem 3.1 of \cite{Bak}. The most difficult part of the proof
of Theorem 5.1 is the proof of item 2. This item is proven in \cite[Theorem 2%
]{BKN}. Item 3 follows from item 2 and an obvious combination of Lemma 2.1
and Theorem 2.1 of \cite{Bak}. Therefore, we do not prove Theorem 5.1 here.

\subsection{The accuracy of the minimizer}

\label{sec:5.2}

We now estimate the accuracy of the minimizer $w_{\min }\in P,$ which is
found in Theorem 5.1. To do this, we recall first that one of main
assumptions of the regularization theory is the \emph{a priori} assumption
of the existence of the exact solution of an ill-posed problem with the
exact, i.e. noiseless data \cite{BK,T}.

Hence, let $c^{\ast }\left( x\right) $ be the exact solution of our CIP with
the exact boundary data $s^{\ast }(x,t),p^{\ast }(x,t)$ in (\ref{2.9}), (\ref%
{2.11}). In this case the function $u(x,t)$ should be replaced with the
function $u^{\ast }(x,t)\in H^{6}$ $\left( Q_{T}\right) $ in (\ref{2.7})-(%
\ref{2.11}), and the function $c\left( x\right) $ in (\ref{2.7}) should be
replaced with $c^{\ast }\left( x\right) ,$ where the function $c^{\ast
}\left( x\right) $ satisfies conditions (\ref{2.1})-(\ref{4.1}). The initial
condition $f\left( x\right) $ in (\ref{2.7}) remains the same. Let 
\begin{equation}
u_{tt}^{\ast }=w^{\ast }\in H_{0}^{4}\left( Q_{T}\right) ,w^{\ast }\mid
_{S_{T}}=q^{\ast },\partial _{\nu }w^{\ast }\mid _{S_{T}}=r^{\ast }.
\label{5.5}
\end{equation}%
We assume that 
\begin{equation}
w^{\ast }\in P^{\ast }=P^{\ast }\left( x_{0},a,b,R,f,q^{\ast },r^{\ast
}\right) ,  \label{5.50}
\end{equation}%
where the set $P^{\ast }$ is obtained from the set $P$ in (\ref{3.13}) via
the replacement of functions $q$ and $r$ with functions $q^{\ast }$ and $%
r^{\ast }$ respectively.

Let a small number $\delta \in \left( 0,1\right) $ be the level of the error
in the boundary data. More precisely, we assume now that there exist
functions $F,F^{\ast }\in H_{0}^{4}\left( Q_{T}\right) $ such that 
\begin{equation}
F\mid _{S_{T}}=q\left( x,t\right) ,\partial _{\nu }F\mid _{S_{T}}=r\left(
x,t\right) ,  \label{5.6}
\end{equation}%
\begin{equation}
F^{\ast }\mid _{S_{T}}=q^{\ast }\left( x,t\right) ,\partial _{\nu }F^{\ast
}\mid _{S_{T}}=r^{\ast }\left( x,t\right) ,  \label{5.7}
\end{equation}%
\begin{equation}
\left\Vert F\right\Vert _{H^{4}\left( Q_{T}\right) },\left\Vert F^{\ast
}\right\Vert _{H^{4}\left( Q_{T}\right) }<R,  \label{5.8}
\end{equation}%
\begin{equation}
\left\Vert F-F^{\ast }\right\Vert _{H^{4}\left( Q_{T}\right) }<\delta <1.
\label{5.9}
\end{equation}%
It is clear from (\ref{2.06}), (\ref{3.3}) and (\ref{3.120}) that 
\begin{equation}
\left\vert \left( \nabla A\left( w\right) ,x-x_{0}\right) \right\vert \leq
B,\forall x\in \overline{\Omega },\forall w\in P.  \label{5.800}
\end{equation}%
Here and below $B=B\left( x_{0},a,b,R,f,Q_{T}\right) >0$ denotes different
positive constants depending only on listed parameters. Choose a small
number $s_{0}>0.$ We impose on the function $w^{\ast }$ a slightly stronger
assumption than the one in the fourth line of (\ref{3.13}). More precisely,
we assume that%
\begin{equation}
\left( \nabla A\left( w^{\ast }\right) ,x-x_{0}\right) \geq s_{0}.
\label{5.901}
\end{equation}

\textbf{Remark 5.1}. \emph{Even though we assume in (\ref{5.9}) that the
noise is less than }$\delta $\emph{\ in the }$H^{4}\left( Q_{T}\right) -$%
\emph{norm, we need this assumption only for the theory. Indeed, in our
computations we introduce a random noise in the boundary data and do not
extend these data inside the time cylinder }$Q_{T}.$\emph{\ Computational
results are successful, see section 6. In other words, computations are less
pessimistic than the theory is. The latter has been always observed in all
previous above cited publications about the convexification, including the
most challenging cases of experimentally collected data \cite%
{Khoa2,KLNSN,KlibSAR2}.}

For any function $w\in P\left( a,b,x_{0},R,f,q,r\right) $ as well as for the
function

$w^{\ast }\in P^{\ast }\left( a,b,x_{0},R,f,q^{\ast },r^{\ast }\right) $
consider functions 
\begin{equation}
v_{\min }=w_{\min }-F,\text{ }v^{\ast }=w^{\ast }-F^{\ast }.  \label{5.9010}
\end{equation}%
It follows from (\ref{2.06}), (\ref{3.50}), (\ref{3.6}), (\ref{3.120}) and (%
\ref{5.9})-(\ref{5.901}) that 
\[
\left( \nabla A\left( v^{\ast }+F\right) ,x-x_{0}\right) =\left( \nabla
A\left( v^{\ast }+F^{\ast }+\left( F-F^{\ast }\right) \right)
,x-x_{0}\right) 
\]%
\begin{equation}
=\left( \nabla A\left( w^{\ast }\right) ,x-x_{0}\right) +V\left( w^{\ast
},F,F^{\ast },x\right) \geq s_{0}-B\delta ,  \label{5.902}
\end{equation}%
where numbers $B$ and $s_{0}$ are defined in (\ref{5.800}) and (\ref{5.901})
respectively. Choose the number $\delta _{1}=\delta _{1}\left(
a,b,R,f,x_{0},s_{0},Q_{T}\right) \in \left( 0,1\right) $ so small that 
\begin{equation}
s_{0}-B\delta \geq 0,\forall \delta \in \left( 0,\delta _{1}\right) ,
\label{5.903}
\end{equation}%
Assume that in (\ref{5.9}) $\delta \in \left( 0,\delta _{1}\right) $ as in (%
\ref{5.903}). Then (\ref{3.50}), (\ref{4.40}), (\ref{3.13}), (\ref{5.50})-(%
\ref{5.903}) and the triangle inequality imply that 
\begin{equation}
v_{\min },v^{\ast }\in P_{0}=P_{0}\left( a,b,x_{0},2R,f,F\right) =\left\{ 
\begin{array}{c}
v\left( x,t\right) :v\in H_{0}^{4,0}\left( Q_{T}\right) , \\ 
\left\Vert v\right\Vert _{H^{4}\left( Q_{T}\right) }\leq 2R, \\ 
a/CR\leq A\left( v+F\right) \leq CRb/a, \\ 
\min_{\overline{\Omega }}\left( \nabla A\left( v+F\right) \left( x\right)
,x-x_{0}\right) \geq 0.%
\end{array}%
\right.  \label{5.90}
\end{equation}

Introduce the functional $I_{\lambda ,\alpha }:P_{0}\rightarrow \mathbb{R}$
as (see (\ref{5.1})):%
\[
I_{\lambda ,\alpha }\left( v\right) =J_{\lambda ,\alpha }\left( v+F\right) 
\]%
\begin{equation}
=\int\limits_{Q_{T}}\left[ A\left( v+F\right) \left( v+F\right) _{tt}-\Delta
\left( v+F\right) \right] ^{2}\varphi _{\lambda }^{2}dxdt+  \label{5.10}
\end{equation}%
\[
+\alpha \left\Vert v+F\right\Vert _{H^{4}\left( Q_{T}\right) }^{2},\text{ }%
v\in P_{0}. 
\]

The following proposition is obvious:

\textbf{Proposition 5.1}. \emph{Let }$\lambda _{1}=\lambda _{1}\left( \eta
,a,b,f,R,P,Q_{T}\right) \geq 1$\emph{\ be the number of Theorem 5.1. The
obvious analog of Theorem 5.1 is valid for functional (\ref{5.10}) for }$%
\lambda \geq \lambda _{2},$ where \emph{\ } 
\begin{equation}
\lambda _{2}=\lambda _{2}\left( \eta ,x_{0},a,b,f,R,P_{0},F^{\ast
},Q_{T}\right) \geq \lambda _{1}\left( \eta ,x_{0},a,b,f,2R,P,Q_{T}\right)
\geq 1.  \label{5.11}
\end{equation}%
\emph{\ }

We have $F^{\ast }$ instead of $F$ in (\ref{5.11}) due to (\ref{5.9}).

\textbf{Theorem 5.2} (the accuracy of the minimizer). \emph{Let the numbers }%
$\eta $\emph{\ and }$T$\emph{\ be the same as in Theorem 4.1. Assume that
the function }$w^{\ast }\in P^{\ast }$ \emph{and conditions (\ref{5.6})-(\ref%
{5.9}), (\ref{5.901}) and (\ref{5.903}) hold. Then:}

\emph{1. In the functional }$I_{\lambda ,\alpha }\left( v\right) $\emph{,
let }$\lambda \geq \lambda _{2},$\emph{\ where the number }$\lambda _{2}$%
\emph{\ is as in (\ref{5.11}). Let }$\alpha \in \left[ 2e^{-\lambda
N},1\right) ,$\emph{\ where the number }$N$\emph{\ is defined in (\ref{4.3}%
). Let }$v_{\min }\in P_{0}$\emph{\ be the minimizer of this functional on
the set }$P_{0}.$\emph{\ Denote }$\overline{w}_{\min }=v_{\min }+F.$\emph{\
Let the function }$c_{\min }\left( x\right) $\emph{\ be the one
reconstructed by the following analog of (\ref{3.9}):}%
\begin{equation}
c_{\min }\left( x\right) =\frac{\Delta f}{\overline{w}_{\min }\left(
x,0\right) }.  \label{5.12}
\end{equation}%
\emph{\ Then the following accuracy estimates hold:}%
\begin{equation}
\left\Vert \overline{w}_{\min }-w^{\ast }\right\Vert _{H^{1}\left(
Q_{T}\right) }\leq C_{3}\left( e^{\lambda \left( M+D^{2}\right) }\delta +%
\sqrt{\alpha }\right) ,  \label{5.13}
\end{equation}%
\begin{equation}
\left\Vert c_{\min }-c^{\ast }\right\Vert _{L_{2}\left( \Omega \right) }\leq
C_{3}\left( e^{\lambda \left( M+D^{2}\right) }\delta +\sqrt{\alpha }\right) ,
\label{5.14}
\end{equation}%
\emph{where numbers }$D^{2}$\emph{\ and }$M$\emph{\ are defined in (\ref{4.2}%
) and (\ref{4.03}) respectively. }

\emph{2. Assume that the number }$\lambda _{2}$ \emph{is so large that }$%
e^{-2\lambda _{2}\left( M+D^{2}\right) }\leq \delta _{1},$\emph{\ where }$%
\delta _{1}$\emph{\ is as in (\ref{5.903}). Let the number }$\delta
_{0}=\delta _{0}\left( \eta ,x_{0},a,b,s_{0},f,R,P_{0},F^{\ast
},Q_{T}\right) \in \left( 0,\delta _{1}\right) \subset \left( 0,1\right) $%
\emph{\ be so small that }%
\begin{equation}
\delta _{0}=e^{-2\lambda _{2}\left( M+D^{2}\right) }.  \label{5.15}
\end{equation}%
\emph{For any }$\delta \in \left( 0,\delta _{0}\right) ,$\emph{\ choose }%
\begin{equation}
\lambda =\lambda \left( \delta \right) =\ln \delta ^{-1/\left( 2\left(
M+D^{2}\right) \right) },\text{ }\alpha =2e^{-\lambda N}.  \label{5.16}
\end{equation}%
\emph{Define the number }$\rho $\emph{\ as} 
\begin{equation}
\rho =\frac{1}{2}\min \left( \frac{N}{M+D^{2}},1\right) \in \left( 0,\frac{1%
}{2}\right] .  \label{5.17}
\end{equation}%
\emph{Then accuracy estimates (\ref{5.13}) and (\ref{5.14}) become} 
\begin{equation}
\left\Vert \overline{w}_{\min }-w^{\ast }\right\Vert _{H^{1}\left(
Q_{T}\right) }\leq C_{3}\delta ^{\rho },  \label{5.18}
\end{equation}%
\begin{equation}
\left\Vert c_{\min }-c^{\ast }\right\Vert _{L_{2}\left( Q_{T}\right) }\leq
C_{3}\delta ^{\rho }.  \label{5.19}
\end{equation}%
\emph{Here and below }$C_{3}=C_{3}\left( \eta
,x_{0},a,b,s_{0},f,R,P_{0},F^{\ast },Q_{T}\right) >0$\emph{\ denotes
different numbers depending only on listed parameters. }

\textbf{Proof}. Recall that by (\ref{5.90}) the function $v^{\ast }\in P_{0}$%
. Hence, we can use Theorem 5.1, (\ref{5.10}) and Proposition 5.1 as:%
\begin{equation}
I_{\lambda ,\alpha }\left( v^{\ast }\right) -I_{\lambda ,\alpha }\left(
v_{\min }\right) -I_{\lambda ,\alpha }\left( v_{\min }\right) \left( v^{\ast
}-v_{\min }\right)  \label{5.20}
\end{equation}%
\[
\geq C_{3}e^{-2\lambda M}\left\Vert v^{\ast }-v_{\min }\right\Vert
_{H^{1}\left( Q_{T}\right) }^{2}. 
\]%
By (\ref{5.4}) $-I_{\lambda ,\alpha }\left( v_{\min }\right) \left( v^{\ast
}-v_{\min }\right) \leq 0.$ Since $-I_{\lambda ,\alpha }\left( v_{\min
}\right) \leq 0$ as well, then (\ref{5.20}) implies that 
\begin{equation}
C_{3}e^{2\lambda M}I_{\lambda ,\alpha }\left( v^{\ast }\right) \geq
\left\Vert v^{\ast }-v_{\min }\right\Vert _{H^{1}\left( Q_{T}\right) }^{2}.
\label{5.21}
\end{equation}%
Since the function $w^{\ast }$ satisfies equation (\ref{3.10}), then (\ref%
{5.1}) implies that $J_{\lambda ,\alpha }\left( w^{\ast }\right) =\alpha
\left\Vert w^{\ast }\right\Vert _{H^{4}\left( Q_{T}\right) }^{2}.$ Next,
using (\ref{4.30}), (\ref{5.9}) and (\ref{5.10}), representing $F=F^{\ast
}+\left( F-F^{\ast }\right) $ and also using (\ref{3.50}), we obtain 
\begin{equation}
I_{\lambda ,\alpha }\left( v^{\ast }\right) =J_{\lambda ,\alpha }\left(
w^{\ast }\right) +\widetilde{I}_{\lambda ,\alpha }\left( v^{\ast }\right)
=\alpha \left\Vert w^{\ast }\right\Vert _{H^{4}\left( Q_{T}\right) }^{2}+%
\widetilde{I}_{\lambda ,\alpha }\left( v^{\ast }\right) ,  \label{5.210}
\end{equation}%
where%
\begin{equation}
\left\vert \widetilde{I}_{\lambda ,\alpha }\left( v^{\ast }\right)
\right\vert \leq C_{3}\delta ^{2}\mathop{\displaystyle \int}%
\limits_{Q_{T}}\varphi _{\lambda }^{2}dxdt\leq C_{3}e^{2\lambda D^{2}}\delta
^{2}.  \label{5.211}
\end{equation}%
It follows from (\ref{5.210}) and (\ref{5.211}) that 
\[
I_{\lambda ,\alpha }\left( v^{\ast }\right) \leq C_{3}\left( e^{2\lambda
D^{2}}\delta ^{2}+\alpha \right) . 
\]%
Combining this with (\ref{5.21}), we obtain (\ref{5.13}).

To obtain (\ref{5.14}), we note that by (\ref{3.9}) and (\ref{5.12}) 
\begin{equation}
c_{\min }\left( x\right) -c^{\ast }\left( x\right) =\Delta f\left( x\right) 
\frac{w^{\ast }\left( x,0\right) -\overline{w}_{\min }\left( x,0\right) }{%
\overline{w}_{\min }\left( x,0\right) w^{\ast }\left( x,0\right) },
\label{5.22}
\end{equation}%
\begin{equation}
\left\Vert u\left( x,0\right) \right\Vert _{L_{2}\left( Q_{T}\right) }\leq
C\left\Vert u\right\Vert _{H^{1}\left( Q_{T}\right) },\forall u\in
H^{1}\left( Q_{T}\right) .  \label{5.23}
\end{equation}%
Thus, (\ref{3.8}), (\ref{5.13}), (\ref{5.22}) and (\ref{5.23}) imply (\ref%
{5.14}). Estimates (\ref{5.18}) and (\ref{5.19}) follows from (\ref{5.13})-(%
\ref{5.17}). $\square $

\subsection{Gradient projection method}

\label{sec:5.3}

Lemma 3.1 implies that the set $P_{0}=P_{0}\left( x_{0},a,b,R,f,F\right) $
defined in (\ref{5.90}) is convex in the space $H_{0}^{4,0}\left(
Q_{T}\right) .$ Therefore, there exists the projection operator $%
Y:H_{0}^{4,0}\left( Q_{T}\right) \rightarrow P_{0}$ of the space $%
H_{0}^{4,0}\left( Q_{T}\right) $ on the set $P_{0}$ \cite[Chapter 10,
section 3.8]{Minoux}. We arrange the gradient projection method of the
minimization of the functional $I_{\lambda ,\alpha }\left( v\right) $ on the
set $P_{0}$ as follows. Let $v_{0}\in P_{0}$ be an arbitrary point. Let the
number $\gamma \in \left( 0,1\right) .$ Then we set 
\begin{equation}
v_{n}=Y\left( v_{n-1}-\gamma I_{\lambda ,\alpha }^{\prime }\left(
v_{n-1}\right) \right) ,\text{ }n=1,2,...  \label{5.24}
\end{equation}%
Note that since by Theorem 5.1 $I_{\lambda ,\alpha }^{\prime }\left(
v_{n-1}\right) \in H_{0}^{4,0}\left( Q_{T}\right) ,$ then all functions $%
v_{n}\in H_{0}^{4,0}\left( Q_{T}\right) .$ Theorem 5.3 follows immediately
from an obvious combination of Theorem 2.1 of \cite{Bak} with Theorems 5.1
and 5.2. The convergence in this theorem is global because the starting
point $v_{0}$ of sequence (\ref{5.24}) is an arbitrary point of the set $%
P_{0}$ and a smallness condition is not imposed on the number $R$ in (\ref%
{5.90}), see Introduction.

\textbf{Theorem 5.3} (global convergence of the gradient projection method (%
\ref{5.24})). \emph{Assume that conditions and notations of Theorem 5.2
hold. Denote }$w_{n}=v_{n}+F,n=0,1,...$\emph{\ Also, let functions }$c_{n,%
\text{comp}}\left( x\right) $\emph{\ be the ones computed by the following
analog of formula (\ref{3.14}):}%
\[
c_{n}\left( x\right) =\frac{\Delta f\left( x\right) }{w_{n}\left( x,0\right) 
}. 
\]%
\emph{\ Then there exists a sufficiently small number }$\gamma \in \left(
0,1\right) $\emph{\ and a number }$\theta =\theta \left( \gamma \right) \in
\left( 0,1\right) $\emph{\ such that the following convergence estimates
hold:} 
\[
\left\Vert w_{n}-\overline{w}_{\min }\right\Vert _{H^{4}\left( Q_{T}\right)
}\leq \theta ^{n}\left\Vert w_{0}-\overline{w}_{\min }\right\Vert
_{H^{4}\left( Q_{T}\right) }, 
\]%
\[
\left\Vert w_{n}-w^{\ast }\right\Vert _{H^{1}\left( Q_{T}\right) }\leq
C_{3}\delta ^{\rho }+\theta ^{n}\left\Vert w_{0}-\overline{w}_{\min
}\right\Vert _{H^{4}\left( Q_{T}\right) }, 
\]%
\[
\left\Vert c_{n}-c^{\ast }\right\Vert _{L_{2}\left( \Omega \right) }\leq
C_{3}\delta ^{\rho }+\theta ^{n}\left\Vert w_{0}-\overline{w}_{\min
}\right\Vert _{H^{4}\left( Q_{T}\right) }. 
\]

\section{Numerical Studies}

\label{sec:6}

\subsection{Numerical implementation}

\label{sec6.1}

\subsubsection{The forward problem (\protect\ref{2.7})-(\protect\ref{2.9})}

\label{6.1.1}

In section 6 $x=\left( x,y,z\right) .$ In our numerical tests, the
computational domain is $\Omega =(0,1)\times (0,1)\times (0,1)$ and $%
x_{0}=(0.5,0.5,-0.5)$ in \eqref{4.1}. To generate the data (\ref{2.11}) for
the CIP, we have solved the forward problem (\ref{2.7})-(\ref{2.9}) by the
finite difference method. We set the total time $T=1/4$. The mesh sizes with
respect to spatial and time variables were the spatial mesh size $h=1/32$ in
all three directions and the mesh size in time is $\tau =1/512$
respectively. Consider two partitions of the interval $[0,1]$ and $[0,1/4]$, 
\begin{eqnarray*}
0 &=&x_{0}<x_{1}<...<x_{32}=1,~x_{p}-x_{p-1}=h, \\
0 &=&y_{0}<y_{1}<...<y_{32}=1,~y_{p}-y_{p-1}=h, \\
0 &=&z_{0}<z_{1}<...<z_{32}=1,~z_{p}-z_{p-1}=h, \\
0 &=&t_{0}<t_{1}<...<t_{128}=1/4,~t_{p}-t_{p-1}=\tau .
\end{eqnarray*}%
Denote $u_{i,j,l}^{k}=u(x_{i},y_{j},z_{l},t_{k})$. We have used an implicit
central difference scheme in time scale and central difference scheme in
space scale. We formulate the implicit discrete scheme as: 
\begin{eqnarray}
c(x_{i},y_{j},z_{l}) &&\frac{u_{i,j,l}^{k+1}+u_{i,j,l}^{k-1}-2u_{i,j,l}^{k}}{%
\tau ^{2}}=\frac{u_{i+1,j,l}^{k^{\ast }}+u_{i-1,j,l}^{k^{\ast
}}-2u_{i,j,l}^{k^{\ast }}}{h^{2}}  \label{6.0} \\
&+&\frac{u_{i,j+1,l}^{k^{\ast }}+u_{i,j-1,l}^{k^{\ast }}-2u_{i,j,l}^{k^{\ast
}}}{h^{2}},\frac{u_{i,j,l+1}^{k^{\ast }}+u_{i,j,l-1}^{k^{\ast
}}-2u_{i,j,l}^{k^{\ast }}}{h^{2}}  \nonumber
\end{eqnarray}%
where $i,j,l$ are the indices with respect to $x,y,z$ and $k$ is the index
respect to time, where $u^{k^{\ast }}=(u^{k-1}+u^{k+1})/2$. The coefficient $%
c(x)$ is assigned as a piecewise constant function to simulate the targets
to be reconstructed, see the next section. As soon as problem {(\ref{2.7})-(%
\ref{2.9}) is solved, we obtain the data for the inverse problem, i.e. the
function }$p\left( x,t\right) $ in {(\ref{2.11}) in the discrete form. In
all our numerical tests, we set in (\ref{2.8}), (\ref{2.9}) }$%
f(x)=x^{2}+y^{2}+z^{2}$ for $x\in \Omega $ and $s(x,t)=f(x)(t+1)$ for $%
(x,t)\in S_{T}$ and use their discrete analogs in computations.

\subsubsection{The inverse problem}

\label{6.1.2}

The required second order $t-$ derivatives of the boundary data in {(\ref%
{3.1})} to solve the inverse problem, see {(\ref{3.12}),} are computed by
the finite differences with respect to $t$ on the boundary $S_{T}$. To solve
the inverse problem, we numerically minimize functional {(\ref{5.1}). To do
this, we write the differential operator in (\ref{5.1}) in the finite
differences and minimize it then with respect to the values of the function }%
${w}${\ at grid points. Thus, we minimize the following functional for }$w^{%
\overline{h}}\in P^{\overline{h}}${\ } 
\begin{equation}
J_{\lambda ,\alpha }^{\overline{h}}\left( w^{\overline{h}}\right)
=\int\limits_{Q_{T}}\left[ A\left( w^{\overline{h}}\right) \left( w^{%
\overline{h}}\right) _{tt}-\Delta \left( w^{\overline{h}}\right) \right]
^{2}\left( \varphi ^{\overline{h}}\right) _{\lambda }^{2}dxdt+\alpha
\left\Vert w^{h}\right\Vert _{H^{2}\left( Q_{T}\right) }^{2},\text{ }
\label{6.1}
\end{equation}%
where the superscript $\overline{h}$ means that derivatives in {(\ref{6.1})
are written via finite differences. }The integration in {(\ref{6.1}) is the
conventional summation over grid points. }To avoid the \textquotedblleft
inverse crime" as well as to decrease the smoothness requirement from $%
H^{4}\left( Q_{T}\right) $ to $H^{2}\left( Q_{T}\right) $, we choose the
mesh step sizes $\widetilde{h}$ in space and $\widetilde{\tau }$ in time
larger than those in the forward problem, $\widetilde{h}=1/16=2h$ in all
three directions, and $\widetilde{\tau }=1/64=8\tau .$

Note that even though the $H^{4}(Q_{T})-$norm was used in the regularization
term of {(\ref{5.1})}, it was used for the theory only. In computations,
however, we use the discrete analog of a simpler $H^{2}(Q_{T})-${norm. This
replacement works well numerically. We explain this by the fact that we use
a rather coarse mesh in (\ref{6.1}), and all norms in finite dimensional
spaces are equivalent. Intuitively, as long as the mesh is rather coarse,
one can still use the equivalence property of norms in finite dimensional
spaces. We have used similar replacements in many of the above cited
publications on the convexification, and they worked well.}

To minimize functional (\ref{6.1}), we apply the conjugate gradient method.
Let $w_{s}^{\overline{h}}$ be the resulting discrete function $w^{\overline{h%
}}$ at the iteration number $s$. Then $w_{s+1}^{\overline{h}}=w_{s}^{%
\overline{h}}+a_{s}d_{s}^{\overline{h}},$ where $a_{s}$ is the iterative
step size and $d_{s}^{\overline{h}}$ is the conjugate gradient direction.
Numbers $a_{s}$ are determined using the line search. The conjugate gradient
directions $d_{s}^{\overline{h}}$ are determined iteratively as follows: Let 
$g_{s}^{\overline{h}}$ be the gradient of the discrete functional $%
J_{\lambda ,\alpha }^{\overline{h}}$ at the point $w_{s}$, i.e. $g_{s}^{%
\overline{h}}=\nabla J_{\lambda ,\alpha }^{\overline{h}}(w_{s}^{\overline{h}%
})$. At the iteration number $0$ we have $d_{0}^{\overline{h}}=-g_{0}^{%
\overline{h}}$ and $\alpha _{0}$ is determined by the line search. Then $%
d_{s+1}^{\overline{h}}=-g_{s+1}^{\overline{h}}+\beta _{s}d_{s}^{\overline{h}%
},$ where the step size $\beta _{s}$ is chosen by Polak-Ribiere-Polyak
formula \cite{Polak} $\beta _{s}=\left[ \left( g_{k+1}^{\overline{h}}\right)
^{T}(g_{k+1}^{\overline{h}}-g_{k}^{\overline{h}})\right] /\Vert g_{k}^{%
\overline{h}}\Vert ^{2}.$ The iterative process stops at the iteration
number $s_{0},$at which $\Vert g_{s_{0}}^{\overline{h}}\Vert _{L_{\infty
}}<\epsilon ,$ where the discrete $L_{\infty }-$norm is taken. In our
numerical experiments, the tolerance number is $\epsilon =10^{-2}$. The
starting point $w_{0}^{\overline{h}}$ for this iterative procedure is the
background data, i.e. the function $w_{0}^{\overline{h}}$ which corresponds
to the solution of forward problem {(\ref{2.7})-(\ref{2.9}) with }$%
c(x)\equiv 1$. Finally, following (\ref{3.14}), we set for computed
functions: $w_{\text{comp}}^{\overline{h}}\left( x,0\right) =w_{s_{0}}^{%
\overline{h}}(x,0),$ 
\[
c_{\text{comp}}^{\overline{h}}(x)=\frac{(\Delta f)(x)}{w_{\text{comp}}^{%
\overline{h}}\left( x,0\right) }. 
\]

\textbf{Remark 6.1}. \emph{We point out that even though the global
convergence of the gradient projection method is claimed in Theorem 5.3, we
have established in our computations that the simpler to implement conjugate
gradient method actually works well. Similar observations took place in all
above cited numerical works of this research group on the convexification.
An explanation of this can be found in \cite[Theorem 4]{KlibSAR2}.}

\subsection{Results}

\label{6.2}

In the tests of this subsection, we demonstrate the efficiency of our
numerical method. To see how our method performs for complicated shapes of
inclusions to be computed, we have chosen in Tests 1-4 complicated shapes of
those inclusions. In all four tests, correct values of the function $c\left(
x\right) $ are: $c\left( x\right) =2$ inside of the inclusion and $c\left(
x\right) =1$ outside of it. Therefore, the correct inclusion/background
contrast is $2:1$ in all cases. Below values of all parameters were chosen
by trial and error. In Tests 1-4, $T=0.25$, $\lambda =1$, $\eta =0.1$. Note
that even though our above theorems require large values of the parameter $%
\lambda ,$ the value $\lambda =1$ works quite well computationally. This
observation is similar to botherations in all above cited previous
computational works on the convexification method. It was established in
these publications that $\lambda \in \left[ 1,3\right] $ works well
computationally. Indeed, it is well known that numerical results are often
less pessimistic than analytical ones.

To test the stability of our method with respect to the random noise in the
data, we introduce a random noise in the boundary data $s(x,t),p(x,t)$ as: 
\begin{equation}
u^{noise}|_{S_{T}}=s(x,t)+\sigma \xi _{x}s(x,t),\text{ }(\partial _{\nu
}u)^{noise}|_{S_{T}}=p(x,t)+\sigma \xi _{x}p(x,t),  \label{6.3}
\end{equation}%
where $\sigma $ is the noise level, and $\xi _{x}$ is the random variable
uniformly distributed in $[-1,1]$. Only points $(x,t)\in S_{T}$ from the
finite difference mesh are taken here. Here, the random variable $\xi _{x}$
depends only on discrete points $x\in \partial \Omega .$ The noise was added
to the boundary data $s(x,t)$ only to solve the inverse problem. But it was
not added when solving forward problem {(\ref{2.7})-(\ref{2.9}) to generate
data (\ref{2.11}) for the inverse problem. }

It is well known that the regularization parameter for an ill-posed problem
should depend on the noise level \cite{BK,T}. In computational practice it
also depends on the range of parameters of this problem. Thus, we take $%
\alpha =0$ for noisy free cases of Tests 1,3. In the cases of noisy data of
Tests 2 and 4 with the 3\% noise level we take $\alpha =0.05.$

\begin{figure}[tbp]
\begin{center}
\begin{tabular}{cc}
\includegraphics[width=4cm]{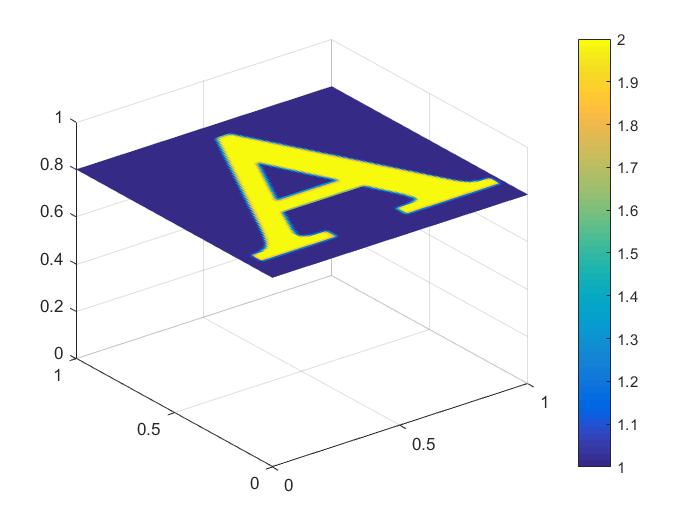} & %
\includegraphics[width=4cm]{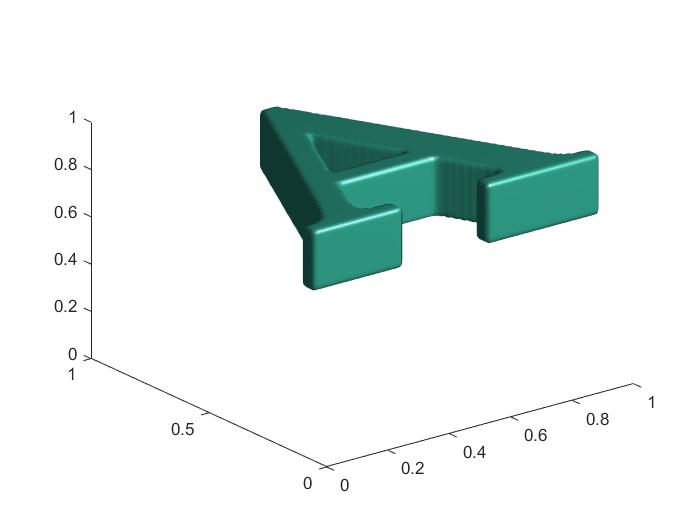} \\ 
(a) Slice image of the true $q$ & (b) 3D image of the true $c$ \\ 
\includegraphics[width=4cm]{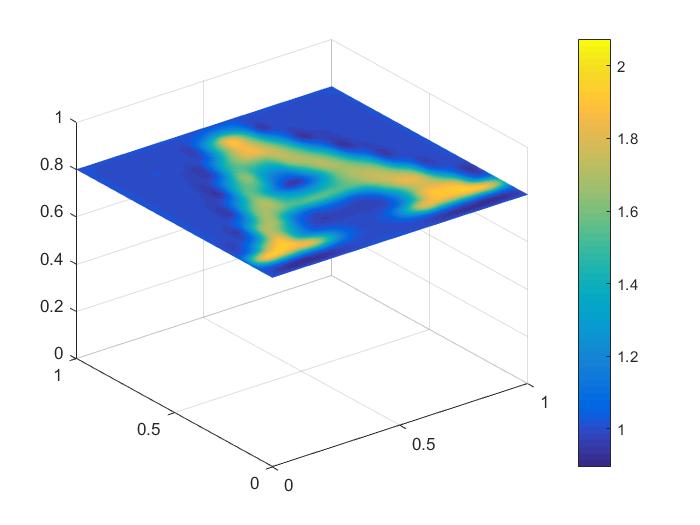} & %
\includegraphics[width=4cm]{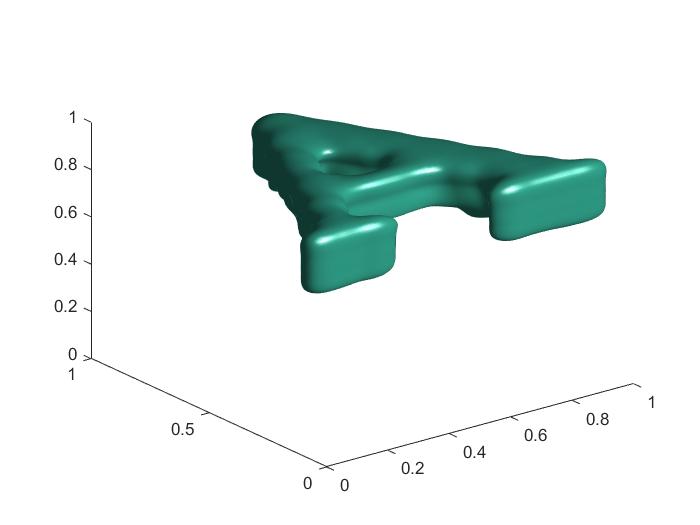} \\ 
(c) Slice image of the recovered $c_{comp}$ & (d) 3D image of the recovered $%
c_{comp}$%
\end{tabular}%
\end{center}
\caption{\emph{Results of Test 1. Noisy free case. The inclusion to be
computed is 'A' shaped with $c =2$ in this shape and $c =1$ outside of it.
(a) and (b) Correct images. (c) and (d) Computed images. One can see from
(c) that $\max c_{\text{comp}}^{\overline{h}}(x)\approx 2,$ which is close
to the true value $c\left( \text{inclusion}\right) =2$. Hence, the
inclusion/background contrast of 2:1 is reconstructed accurately. Comparing
(b) and (d), one can see that the shape of the inclusion is also accurately
reconstructed.}}
\label{example-A}
\end{figure}

\textbf{Test 1}. Noisy free case. We test the reconstruction by our method
of the case when the shape of our inclusion is the same as the shape of the
letter `$A$'. The true function $c(x)$ is depicted on Figures \ref{example-A}
(a), (b). We set $c=2$ inside of this inclusion and $c=1$ outside of it.
Thus, the inclusion/background contrast is 2:1. See Figures \ref{example-A}
for the reconstruction results.

\begin{figure}[tbp]
\begin{center}
\begin{tabular}{cc}
\includegraphics[width=4cm]{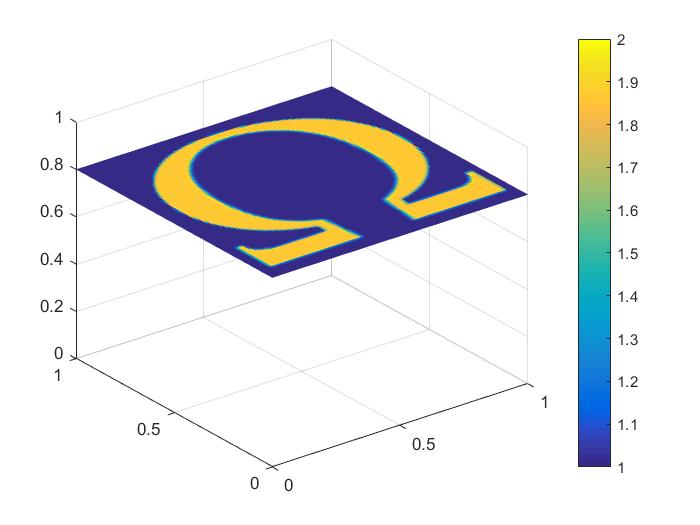} & %
\includegraphics[width=4cm]{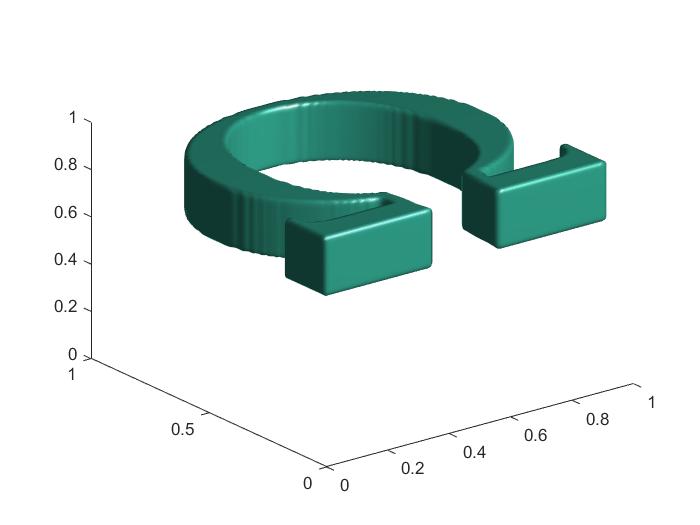} \\ 
(a) Slice image of the true $q$ & (b) 3D image of the true $c$ \\ 
\includegraphics[width=4cm]{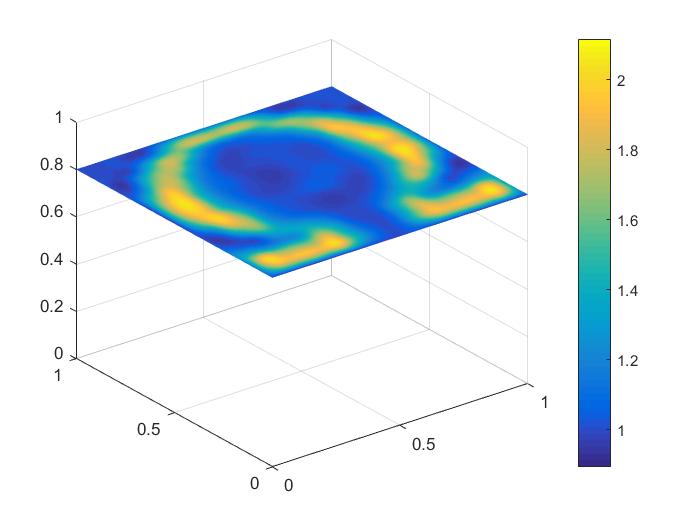}
& %
\includegraphics[width=4cm]{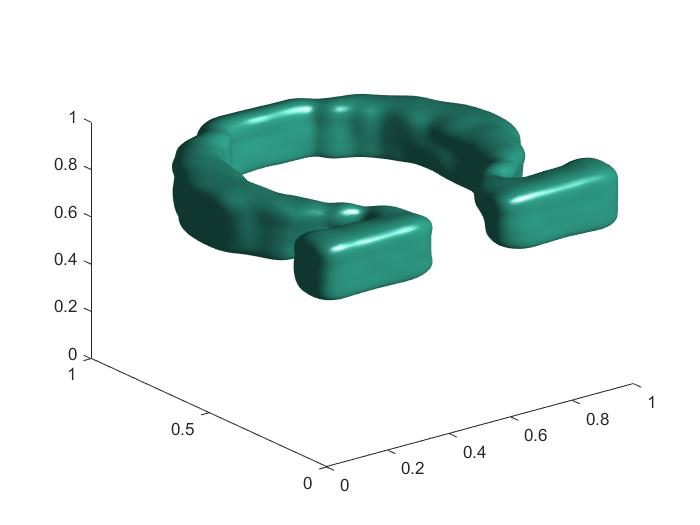}
\\ 
(c) Slice image of the recovered $c_{comp}$ & (d) 3D image of the recovered $%
c_{comp}$%
\end{tabular}%
\end{center}
\caption{\emph{Results of Test 2. Noisy data with $\protect\sigma =0.03$ in (%
\protect\ref{6.3}), i.e. 3\% noise. The inclusion to be computed is '$\Omega$%
' shaped with $c =2$ in this shape and $c =1$ outside of it. (a) and (b)
Correct images. (c) and (d) Computed images. In this example, we add 3$\%$
noise. One can see from (c) that $\max c_{\text{comp}}^{\overline{h}%
}(x)\approx 2,$ which is close to the true value $c\left( \text{inclusion}%
\right) =2$. Hence, the inclusion/background contrast of 2:1 is
reconstructed accurately. Comparing (b) and (d), one can see that the shape
of the inclusion is also accurately reconstructed.}}
\label{example-Omega}
\end{figure}

\textbf{Test 2}. We now use the noisy boundary data (\ref{6.3}) with $\sigma
=0.03,$ which means 3\% of random noise. We test the reconstruction by our
method of the case when the shape of our inclusion is the same as the shape
of the letter `$\Omega $'. The function $c(x)$ is depicted on Figures \ref%
{example-Omega} (a) and (b). See Figures \ref{example-Omega} for the
reconstruction results.

\begin{figure}[tbp]
\begin{center}
\begin{tabular}{cc}
\includegraphics[width=4cm]{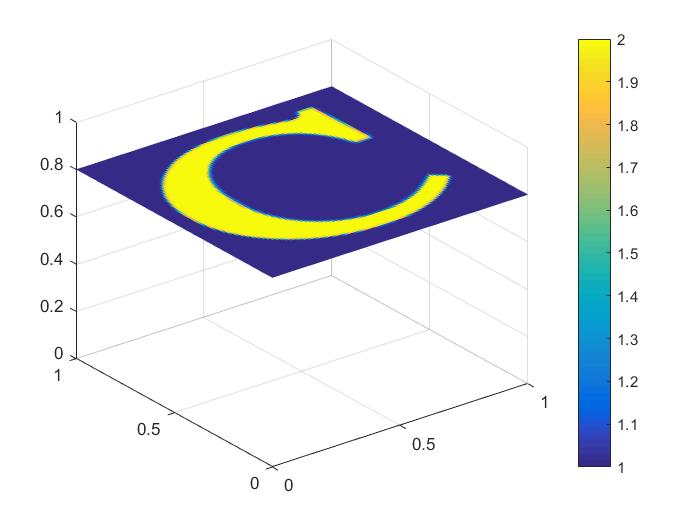} & %
\includegraphics[width=4cm]{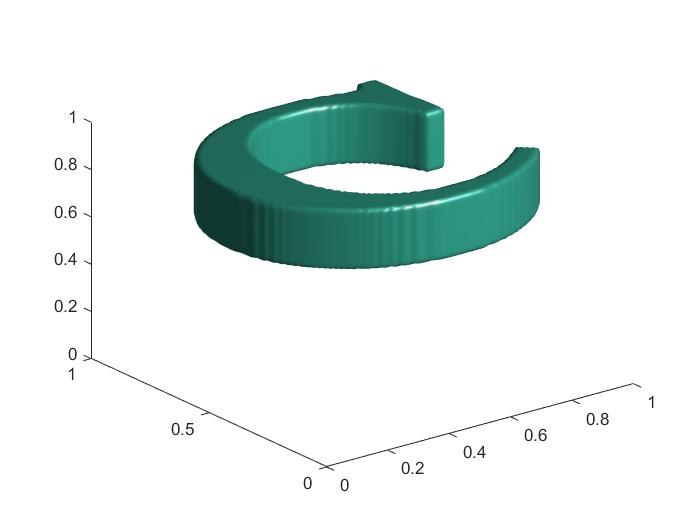} \\ 
(a) Slice image of the true $q$ & (b) 3D image of the true $c$ \\ 
\includegraphics[width=4cm]{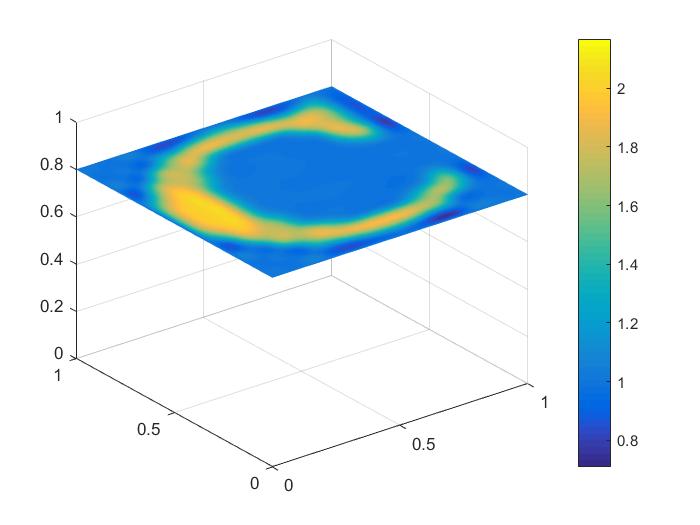} & %
\includegraphics[width=4cm]{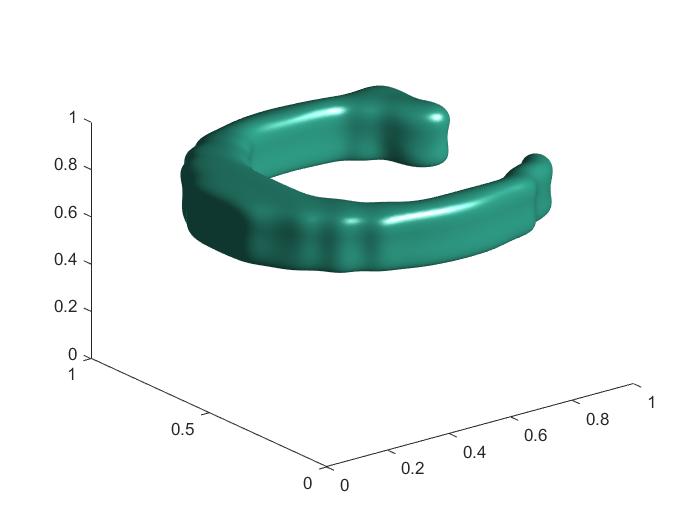} \\ 
(c) Slice image of the recovered $c_{comp}$ & (d) 3D image of the recovered $%
c_{comp}$%
\end{tabular}%
\end{center}
\caption{\emph{Results of Test 3. Noisy free case. The inclusion to be
computed is 'C' shaped with $c =2$ in this shape and $c =1$ outside of it.
(c) and (d) Computed images. One can see from (c) that $\max c_{\text{comp}%
}^{\overline{h}}(x)\approx 2,$ which is close to the true value $c\left( 
\text{inclusion}\right) =2$. Hence, the inclusion/background contrast of 2:1
is reconstructed accurately. Comparing (b) and (d), one can see that the
shape of the inclusion is also accurately reconstructed.}}
\label{example-C}
\end{figure}

\textbf{Test 3}. Noisy free case. We test the reconstruction by our method
of the case when the shape of our inclusion is the same as the shape of the
letter `$C$'. The function $c(x)$ is depicted on Figures \ref{example-C} (a)
and (b). See Figures \ref{example-C} for the reconstruction results.

\begin{figure}[tbp]
\begin{center}
\begin{tabular}{cc}
\includegraphics[width=4cm]{Figures/C-up-0-8-slice-true.jpg} & %
\includegraphics[width=4cm]{Figures/C-up-0-8-sur-true.jpg} \\ 
(a) Slice image of the true $q$ & (b) 3D image of the true $c$ \\ 
\includegraphics[width=4cm]{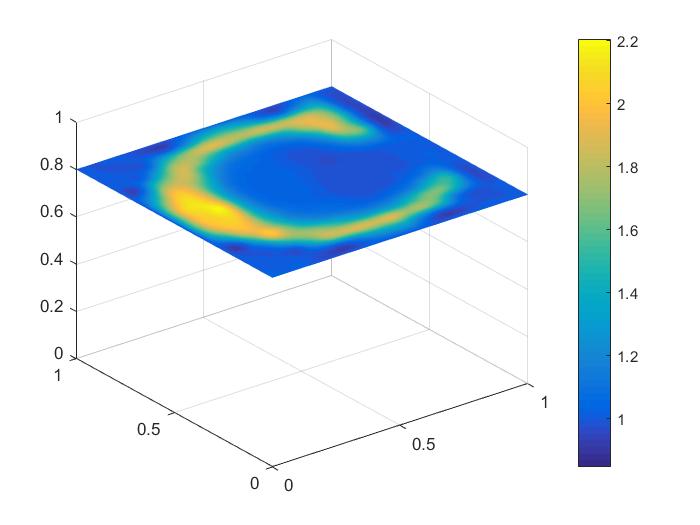} & %
\includegraphics[width=4cm]{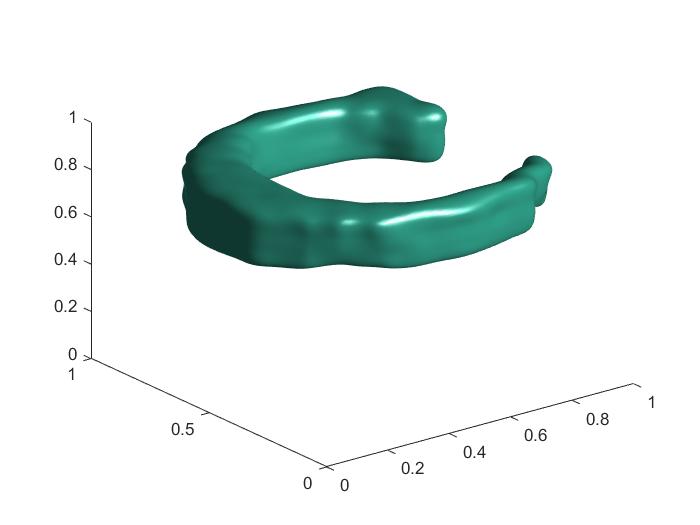} \\ 
(c) Slice image of the recovered $c_{comp}$ & (d) 3D image of the recovered $%
c_{comp}$%
\end{tabular}%
\end{center}
\caption{\emph{Results of Test 4. Noisy data with $\protect\sigma =0.03$ in (%
\protect\ref{6.3}), i.e. 3\% noise. The inclusion to be computed is 'C'
shaped with $c =2$ in this shape and $c =1$ outside of it. (c) and (d)
Computed images. In this example, we add 3$\%$ noise. One can see from (c)
that $\max c_{\text{comp}}^{\overline{h}}(x)\approx 2,$ which is close to
the true value $c\left( \text{inclusion}\right) =2$. Hence, the
inclusion/background contrast of 2:1 is reconstructed accurately. Comparing
(b) and (d), one can see that the shape of the inclusion is also accurately
reconstructed.}}
\label{example-C-noise}
\end{figure}

\textbf{Test 4}. Just as in Test 2, we use the noisy data (\ref{6.3}) with
uniform distributed random variable on $[-1,1]$ of the 3\% level, i.e. $%
\sigma =0.03.$ The shape of our inclusion is the same as the shape of the
letter `$C$'. The function $c(x)$ is depicted on Figures \ref%
{example-C-noise} (a) and (b). See Figures \ref{example-C-noise} for the
reconstruction results.

On Figures 6.1-6.4, each computed function $c_{\text{comp}}^{\overline{h}}$
is depicted by a 2-D slice and by the surface plot. We have not applied any
postprocessing procedure to our computed images.

\end{document}